\documentclass[10pt]{amsart}
\usepackage{amsopn}
\usepackage{amsmath, amsthm, amssymb, amscd, amsfonts, xcolor}
\usepackage[latin1]{inputenc}

\usepackage{xcolor}
\usepackage{enumitem}
\usepackage{multirow}

\input xy
\xyoption{all}

\usepackage[OT2,OT1]{fontenc}
\newcommand\cyr{%
\renewcommand\rmdefault{wncyr}%
\renewcommand\sfdefault{wncyss}%
\renewcommand\encodingdefault{OT2}%
\normalfont
\selectfont}
\DeclareTextFontCommand{\textcyr}{\cyr}

\usepackage[pdftex]{hyperref}

\theoremstyle{plain}

\newtheorem{teo}{Theorem}[section]
\newtheorem{cor}[teo]{Corollary}

\newtheorem{prop}[teo]{Proposition}
\newtheorem{lema}[teo]{Lemma}

\theoremstyle{definition}

\numberwithin{equation}{section}

\newcommand\Ad{\operatorname{Ad}}
\newcommand\ad{\operatorname{ad}}

\newcommand\Aut{\operatorname{Aut}}

\newcommand\rk{\operatorname{rk}}

\newcommand\Exp{\operatorname{Exp}}

\newcommand\codim{\operatorname{codim}}

\newcommand\so{\mathfrak{so}}
\newcommand\tr{\mathfrak{tr}}
\renewcommand\gg{\mathfrak{g}}
\newcommand\ga{\mathfrak{a}}
\newcommand\gb{\mathfrak{b}}
\newcommand\gh{\mathfrak{h}}
\newcommand\gk{\mathfrak{k}}
\newcommand\gl{\mathfrak{l}}
\newcommand\gm{\mathfrak{m}}
\newcommand\gn{\mathfrak{n}}
\newcommand\gf{\mathfrak{f}}
\newcommand\gp{\mathfrak{p}}
\newcommand\gq{\mathfrak{q}}
\newcommand\gz{\mathfrak{z}}
\newcommand\gsl{\mathfrak{sl}}
\newcommand\ggl{\mathfrak{gl}}
\newcommand\gso{\mathfrak{so}}
\newcommand\bbc{\mathbb{C}}
\newcommand\bbv{\mathbb{V}}

\newcommand\bbw{\mathbb{W}}
\newcommand\bbr{\mathbb{R}}
\newcommand\bbh{\mathbb{H}}
\newcommand\bbo{\mathbb{O}}

\makeatletter\renewcommand\theenumi{\@roman\c@enumi}\makeatother

\title{Maximal totally geodesic submanifolds\\ and index of symmetric spaces}

\author{J\"{u}rgen Berndt}
\address{King's College London, Department of Mathematics, London WC2R 2LS, United Kingdom}
\email{jurgen.berndt@kcl.ac.uk}
\thanks{}

\author{Carlos Olmos}
\address{Facultad de Matem\'atica, Astronom\'ia y F\'isica, Universidad Nacional de C\'ordoba, 
Ciudad Universitaria, 5000 C\'ordoba, Argentina}
\email{olmos@famaf.unc.edu.ar}
\thanks{This research was supported by Famaf-UNC and CIEM-Conicet. The article was written while the first author visited the University of California, Irvine. He would like to thank Professor Chuu-Lian Terng and the University for their kind support and hospitality during the visit.}


\begin {document}

\begin{abstract}
Let $M$ be an irreducible Riemannian symmetric space. The index $i(M)$ of $M$ is the minimal codimension of a totally geodesic submanifold of $M$. In \cite{BO} we proved that $i(M)$ is bounded from below by the rank $\rk(M)$ of $M$, that is, $\rk(M) \leq i(M)$. In this paper we classify all irreducible Riemannian symmetric spaces $M$ for which the equality holds, that is, $\rk(M) = i(M)$. In this context we also obtain an explicit classification of all non-semisimple maximal totally geodesic submanifolds in irreducible Riemannian symmetric spaces of noncompact type and show that they are closely related to irreducible symmetric R-spaces. We also determine the index of some symmetric spaces and classify the irreducible Riemannian symmetric spaces of noncompact type with $i(M) \in \{4,5,6\}$.
\end{abstract}

\maketitle 

\section {Introduction}

Let $M$ be a connected Riemannian manifold and denote by ${\mathcal S}$ the set of all connected totally geodesic submanifolds $\Sigma$ of $M$ with $\dim(\Sigma) < \dim(M)$. The index $i(M)$ of $M$ is defined by
\[
i(M) = \min\{ \dim(M) - \dim(\Sigma) \mid \Sigma \in {\mathcal S}\} = \min\{ \codim(\Sigma) \mid \Sigma \in {\mathcal S}\}.
\]
This notion was introduced by Onishchik in \cite{On} who also classified the irreducible simply connected Riemannian symmetric spaces $M$  with $i(M) \leq 2$. 

In \cite{BO} we investigated $i(M)$ for irreducible Riemannian symmetric spaces $M$. We proved that the rank $\rk(M)$ of $M$ is always less or equal than the index of $M$ and classified all irreducible Riemannian symmetric spaces $M$ with $i(M) \leq 3$. The motivation for this paper was to understand better the equality case $\rk(M) = i(M)$. The main result of this paper is the classification of all irreducible Riemannian symmetric spaces $M$ with $\rk(M) = i(M)$.

\begin{teo}\label{main}
Let $M$ be an irreducible Riemannian symmetric space of noncompact type. The equality $\rk(M) = i(M)$ holds if and only if $M$ is isometric to one of the following symmetric spaces:
\begin{itemize}[leftmargin=.3in]
\item[\rm{(i)}] $SL_{r+1}({\mathbb R})/SO_{r+1}$, $r \geq 1$;
\item[\rm{(ii)}] $SO^o_{r,r+k}/SO_rSO_{r+k}$, $r \geq 1$, $k \geq 0$, $(r,k) \notin \{(1,0),(2,0)\}$.
\end{itemize}
\end{teo}

\noindent 
Duality between Riemannian symmetric spaces of noncompact type and of compact type preserves totally geodesic submanifolds, and if $M$ is an irreducible Riemannian symmetric space of compact type and $\hat{M}$ is its Riemannian universal covering space (which is also a Riemannian symmetric space of compact type), then $i(M) = i(\hat{M})$. Therefore Theorem \ref{main} leads, via duality and covering maps, to the classification of irreducible Riemannian symmetric spaces of compact type with $\rk(M) = i(M)$.

In order to compute the index explicitly we need to have a good understanding of maximal totally geodesic submanifolds. Every maximal totally geodesic submanifold $\Sigma$ in an irreducible Riemannian symmetric space $M$ of noncompact type is either semisimple or non-semisimple. As part of our investigation we obtain an explicit classification for the non-semisimple case and a conceptual characterization of such submanifolds in terms of symmetric R-spaces. Denote by $r$ the rank of $M$ and write $M = G/K$, where $G$ is the connected identity component of the isometry group $I(M)$ of $M$ and $K = G_p$ is the isotropy group of $G$ at $p \in M$. Consider a set of simple roots $\Lambda = \{\alpha_1,\ldots,\alpha_r\}$ of a restricted root space decomposition of the Lie algebra $\gg$ of $G$ and denote by $\delta = \delta_1 \alpha_1 + \ldots + \delta_r \alpha_r$ the highest root. Let $\gq_i$ be the parabolic subalgebra of $\gg$ which is determined by the root subsystem $\Phi_i = \Lambda \setminus \{\alpha_i\}$ and consider the Chevalley decomposition $\gq_i = \gl_i \oplus \gn_i$ of $\gq_i$ into a reductive subalgebra $\gl_i$ and a nilpotent subalgebra $\gn_i$. Let $L_i$ be the connected closed subgroup of $G$ with Lie algebra $\gl_i$ and denote by $F_i$ the orbit of $L_i$ containing $p$. Then $F_i$ is a non-semisimple totally geodesic submanifold of $M$ which decomposes into $F_i = \bbr \times B_i$, where $B_i$ is a semisimple Riemannian symmetric space of noncompact type. The classification and characterization of non-semisimple maximal totally geodesic submanifolds in $M$ is as follows:

\begin{teo} \label{main2}
Let $M = G/K$ be an irreducible Riemannian symmetric space of noncompact type and let $\Sigma$ be a non-semisimple connected complete totally geodesic submanifold of $M$. Then the following statements are equivalent:
\begin{itemize}[leftmargin=.3in]
\item[\rm{(i)}] $\Sigma$ is a maximal totally geodesic submanifold of $M$;
\item[\rm{(ii)}] $\Sigma$ is isometrically congruent to $F_i = \bbr \times B_i$ and $\delta_i = 1$;
\item[\rm{(iii)}] The normal space $\nu_p\Sigma$ of $\Sigma$ at $p$ is the tangent space of a symmetric R-space in $T_pM$;
\item[\rm{(iv)}] The pair $(M,\Sigma)$ is as in Table {\rm \ref{totgeodnsstable}}.
\end{itemize}
\end{teo}

An R-space is a real flag manifold and a symmetric R-space is a real flag manifold which is also a symmetric space. R-spaces are projective varieties and symmetric R-spaces were classified and investigated by Kobayashi and Nagano in \cite{KN}. They arise as certain orbits of the isotropy representation of semisimple Riemannian symmetric spaces.

This paper is organized as follows. In Section \ref{preliminaries} we summarize basic material about Riemannian symmetric spaces of noncompact type, their restricted root space decompositions and Dynkin diagrams, parabolic subalgebras, and their boundary components with respect to the maximal Satake compactification.

In Section \ref{reflective} we obtain some sufficient criteria for totally geodesic submanifolds in Riemannian symmetric spaces of noncompact type to be reflective. As is well-known, totally geodesic submanifold are in one-to-one correspondence with Lie triple system. If the orthogonal complement of a Lie triple system is also a Lie triple system, then the Lie triple system and the corresponding totally geodesic submanifold are said to be reflective. Geometrically, reflective submanifolds arise as connected components of fixed point sets of isometric involutions. Reflective submanifolds in irreducible simply connected Riemannian symmetric spaces of compact type were classified by Leung in \cite{L1} and \cite{L2}. The concept of reflectivity turns out to be very useful in our context. One of our main criteria is Proposition \ref{ref3} which states that if the kernel of the slice representation of a semisimple totally geodesic submanifold $\Sigma$ in an irreducible Riemannian symmetric space of noncompact type has positive dimension, then $\Sigma$ is reflective. This criterion then provides a lower bound for the codimension of $\Sigma$ which we will use in index calculations.

In Section \ref{nonsemisimple} we will prove Theorem \ref{main2}. The first step is to show that any non-semisimple maximal totally geodesic submanifold in $M$ is congruent to one of the orbits $F_i$ introduced above. The coefficient $\delta_i$ of $\alpha_i$ in the highest root $\delta$ then plays a crucial role for the next step. If $\delta_i \geq 2$, we construct explicitly a larger Lie triple system containing the Lie triple system corresponding to $F_i$. The situation for $\delta_1 = 1$ is much more involved. With delicate arguments using Killing fields, Jacobi fields, reflections and transvections we can show that $F_i$ is maximal when $\delta_i = 1$. As an application of Theorem \ref{main2} we obtain that every maximal totally geodesic submanifold of an irreducible Riemannian symmetric space of noncompact type whose root system is of type $(BC_r)$, $(E_8)$, $(F_4)$ or $(G_2)$ must be semisimple. Another application states that every non-semisimple maximal totally geodesic submanifold of an irreducible Riemannian symmetric space of noncompact type must be reflective. As a third application we obtain that the index of $SL_{r+1}(\bbr)/SO_{r+1}$ is equal to its rank $r$.

In Section \ref{ex} we prove that the two classes of symmetric spaces listed in Theorem \ref{main} satisfy the equality $\rk(M) = i(M)$. For this we explicitly construct totally geodesic submanifolds $\Sigma$ of $M$ with $\codim(\Sigma) = \rk(M)$ using standard algebraic theory of symmetric spaces. 

In Section \ref{proof} we prove Theorem \ref{main}. A crucial step is Proposition \ref{boundary_reduction} which states that if $M$ satisfies the equality $\rk(M) = i(M)$, then every irreducible boundary component $B$ of the maximal Satake compactification of $M$ satisfies $\rk(B) = i(B)$. As an application we obtain that with the possible exception of $E_6^6/Sp_4$, $E_7^7/SU_8$ and $E_8^8/SO_{16}$ there are no other irreducible Riemannian symmetric spaces $M$ of noncompact type with $\rk(M) = i(M)$ than those discussed in Section \ref{ex}. The exceptional symmetric space $E_6^6/Sp_4$ has the interesting property that each of its irreducible boundary components $B$ satisfies $\rk(B) = i(B)$. In order to come to a conclusion for this exceptional symmetric space we developed the criteria about reflective submanifolds in Section \ref{reflective}. Using these criteria we can show that $E_6^6/Sp_4$ does not satisfy the equality $\rk(M) = i(M)$. Since $E_6^6/Sp_4$ arises as a boundary component of $E_7^7/SU_8$ and of $E_8^8/SO_{16}$ we can then conclude that these two symmetric spaces do not satisfy the equality $\rk(M) = i(M)$ either. 

In Section \ref{applications} we apply some of the results in Sections \ref{reflective} and \ref{nonsemisimple} to calculate explicitly the index of some other symmetric spaces. We also classify the irreducible Riemannian symmetric spaces of noncompact type with $i(M) \in \{4,5,6\}$. 

\section{Riemannian symmetric spaces of noncompact type} \label{preliminaries}

We assume that the reader is familiar with the general theory of Riemannian symmetric spaces as in \cite{H} and summarize below some basic facts and notations which are used in this paper.

Let $M = G/K$ be an irreducible Riemannian symmetric space of noncompact type, where $G = I^o(M)$ is the connected component of the isometry group $I(M)$ of $M$ containing the identity transformation, $p \in M$ and $K = G_p$ is the isotropy group of $G$ at $p$. Then $G$ is a noncompact real semisimple Lie group and $K$ is a maximal compact subgroup of $G$. Let $\gg = \gk \oplus \gp$ be the corresponding Cartan decomposition of $\gg$ and denote by $\theta$ the corresponding Cartan involution on $\gg$. Let $B$ be the Killing form of $\gg$. Then $\langle X,Y \rangle = -B(X,\theta Y)$ is a positive definite inner product on $\gg$. The vector space $\gp$ can be identified canonically with the tangent space $T_pM$ of $M$ a $p$. Since the Riemannian metric on $M$ is unique up to homothety, we can assume that the Riemannian metric on $M$ coincides with the $G$-invariant Riemannian metric  induced by $\langle \cdot , \cdot \rangle$.

We denote by $r = \rk(M)$ the rank of $M$. Let $\ga$ be a maximal abelian subspace of $\gp$ and denote by $\ga^\ast$ the dual space of $\ga$. Note that $\dim(\ga) = r$. For each $\alpha \in \ga^\ast$ we define $\gg_{\alpha} = \{X \in \gg \mid [H,X] = \alpha(H)X\ {\rm for\ all\ }H \in \ga\}$. If $\alpha \neq 0$ and $\gg_\alpha \neq \{0\}$, then $\alpha$ is a restricted root and $\gg_\alpha$ a restricted root space of $\gg$ with respect to $\ga$. The positive integer $m_\alpha = \dim(\gg_\alpha)$ is called the multiplicity of the root $\alpha$. We denote by $\Psi$ the set of restricted roots with respect to $\ga$. The direct sum decomposition
\[
\gg = \gg_0 \oplus \left(\bigoplus_{\alpha \in \Psi} \gg_{\alpha}\right)
\]
is the restricted root space decomposition of $\gg$ with respect to $\ga$. The eigenspace $\gg_0$ decomposes into $\gg_0 = \gk_0 \oplus \ga$, where $\gk_0 = Z_{\gk}(\ga)$ is the centralizer of $\ga$ in $\gk$. 

Let $\{\alpha_1,\ldots,\alpha_r\} = \Lambda \subset \Psi$ be a set of simple roots of $\Psi$. We denote by $H^1,\ldots,H^r \in \ga$ the dual basis of $\alpha_1,\ldots,\alpha_r \in \ga^*$ defined by $\alpha_i(H^j) = \delta_{ij}$ for all $i,j \in \{1,\ldots,r\}$, where $\delta_{ij} = 0$ for $i \neq j$ and $\delta_{ij} = 1$ for $i = j$. Riemannian symmetric spaces of noncompact type are uniquely determined by the Dynkin diagram of their restricted root system together with the multiplicities of the simple roots. In Table \ref{dynkin} we list the Dynkin diagrams and root multiplicities for all irreducible Riemannian symmetric spaces of noncompact type.

\begin{table}
\caption{Dynkin diagrams and root multiplicities for irreducible Riemannian symmetric spaces $M$ of noncompact type} 
\label{dynkin} 
{\footnotesize\begin{tabular}{ | p{4.5cm}  p{2.6cm}  p{2.4cm}  p{1.6cm} | }
\hline \rule{0pt}{4mm}
Dynkin diagram & $M$ & Multiplicities & Comments\\[1mm]
\hline \rule{0pt}{4mm}
\multirow{5}{*}{\tiny$
\xy
\POS (0,0) *\cir<2pt>{} ="a",
(7,0) *\cir<2pt>{}="b",
(14,0) *\cir<2pt>{}="c",
(21,0) *\cir<2pt>{}="d",
(0,-3) *{\alpha_1},
(7,-3) *{\alpha_2},
(14,-3) *{\alpha_{r-1}},
(21,-3) *{\alpha_r},
\ar @{-} "a";"b",
\ar @{.} "b";"c",
\ar @{-} "c";"d",
\endxy
$} 
& $SO^o_{1,1+k}/SO_{1+k}$ & $k$  & $k \geq 1$\\
& $SL_{r+1}({\mathbb R})/SO_{r+1}$ & $1,1,\ldots,1,1$ & $r \geq 2$\\
& $SL_{r+1}({\mathbb C})/SU_{r+1}$ & $2,2,\ldots,2,2$ & $r \geq 2$\\
& $SU^*_{2r+2}/Sp_{r+1}$ & $4,4,\ldots,4,4$ &  $r \geq 2$\\
& $E_6^{-26}/F_4$ & $8,8$ & \\[1mm]
\hline \rule{0pt}{4mm}
\multirow{2}{*}{\tiny$\xy
\POS (0,0) *\cir<2pt>{} ="a",
(7,0) *\cir<2pt>{}="b",
(14,0) *\cir<2pt>{}="c",
(21,0) *\cir<2pt>{}="d",
(28,0) *\cir<2pt>{}="e",
(0,-3) *{\alpha_1},
(7,-3) *{\alpha_2},
(14,-3) *{\alpha_{r-2}},
(21,-3) *{\alpha_{r-1}},
(28,-3) *{\alpha_r},
\ar @{-} "a";"b",
\ar @{.} "b";"c",
\ar @{-} "c";"d",
\ar @2{->} "d";"e"
\endxy$}
& $SO^o_{r,r+k}/SO_{r}SO_{r+k}$ & $1,1,\ldots,1,1,k$ &  $r \geq 2, k \geq 1$\\
& $SO_{2r+1}({\mathbb C})/SO_{2r+1}$ & $2,2,\ldots,2,2,2$ & $r \geq 2$\\[1mm]
\hline \rule{0pt}{4mm}
\multirow{6}{*}{\tiny$\xy
\POS (0,0) *\cir<2pt>{} ="a",
(7,0) *\cir<2pt>{}="b",
(14,0) *\cir<2pt>{}="c",
(21,0) *\cir<2pt>{}="d",
(28,0) *\cir<2pt>{}="e",
(0,-3) *{\alpha_1},
(7,-3) *{\alpha_2},
(14,-3) *{\alpha_{r-2}},
(21,-3) *{\alpha_{r-1}},
(28,-3) *{\alpha_r},
\ar @{-} "a";"b",
\ar @{.} "b";"c",
\ar @{-} "c";"d",
\ar @2{<-} "d";"e"
\endxy$}
& $Sp_r({\mathbb R})/U_r$ & $1,1,\ldots,1,1,1$ & $r \geq 3$\\
& $SU_{r,r}/S(U_rU_r)$ &  $2,2,\ldots,2,2,1$ & $r \geq 3$\\
& $Sp_r({\mathbb C})/Sp_r$ & $2,2,\ldots,2,2,2$ &  $r \geq 3$\\
& $SO^*_{4r}/U_{2r}$ & $4,4,\ldots,4,4,1$ &  $r \geq 3$\\
& $Sp_{r,r}/Sp_rSp_r$ &  $4,4,\ldots,4,4,3$ & $r \geq 2$\\
& $E_7^{-25}/E_6U_1$ & $8,8,1$ & \\[1mm]
\hline \rule{0pt}{4mm}
\multirow{2}{*}{\tiny$\xy
\POS (0,0) *\cir<2pt>{} ="a",
(7,0) *\cir<2pt>{}="b",
(14,0) *\cir<2pt>{}="c",
(21,0) *\cir<2pt>{}="d",
(28,2) *\cir<2pt>{}="e",
(28,-2) *\cir<2pt>{}="f",
(0,-3) *{\alpha_1},
(7,-3) *{\alpha_2},
(14,-3) *{\alpha_{r-3}},
(21,-3) *{\alpha_{r-2}},
(32,3) *{\alpha_{r-1}},
(31,-3) *{\alpha_r},
\ar @{-} "a";"b",
\ar @{.} "b";"c",
\ar @{-} "c";"d",
\ar @{-} "d";"e",
\ar @{-} "d";"f"
\endxy$} & $SO^o_{r,r}/SO_{r}SO_{r}$ & $1,1,\ldots,1,1,1,1$ & $r \geq 4$\\
& $SO_{2r}({\mathbb C})/SO_{2r}$ &  $2,2,\ldots,2,2,2,2$ & $r \geq 4$\\[1mm]
\hline \rule{0pt}{4mm}
\multirow{5}{*}{\tiny$\xy
\POS (0,0) *\cir<2pt>{} ="a",
(7,0) *\cir<2pt>{}="b",
(14,0) *\cir<2pt>{}="c",
(21,0) *\cir<2pt>{}="d",
(30,0) *\cir<2pt>{},
(30,0) *\cir<4pt>{}="e",
(0,-3) *{\alpha_1},
(7,-3) *{\alpha_2},
(14,-3) *{\alpha_{r-2}},
(21,-3) *{\alpha_{r-1}},
(30,-3) *{(\alpha_r,2\alpha_r)},
\ar @{-} "a";"b",
\ar @{.} "b";"c",
\ar @{-} "c";"d",
\ar @2{<->} "d";"e"
\endxy$} & $SU_{r,r+k}/S(U_rU_{r+k})$ & $2,2,\ldots,2,2,(2k,1)$ & $r \geq 1, k \geq 1$\\
& $Sp_{r,r+k}/Sp_rSp_{r+k}$ & $4,4,\ldots,4,4,(4k,3)$ & $r \geq 1, k \geq 1$\\
& $SO^*_{4r+2}/U_{2r+1}$ & $4,4,\ldots,4,4,(4,1)$ &  $r \geq 2$\\
& $E_6^{-14}/Spin_{10}U_1$ & $6,(8,1)$ &  \\
& $F_4^{-20}/Spin_9$ & $(8,7)$ &  \\[1mm]
\hline \rule{0pt}{4mm}
\multirow{2}{*}{\tiny$\xy
\POS (0,0) *\cir<2pt>{} ="a",
(14,5) *\cir<2pt>{} = "b",
(7,0) *\cir<2pt>{}="c",
(14,0) *\cir<2pt>{}="d",
(21,0) *\cir<2pt>{}="e",
(28,0) *\cir<2pt>{}="f",
(0,2) *{\alpha_1},
(11,5) *{\alpha_2},
(7,2) *{\alpha_3},
(16,2) *{\alpha_4},
(21,2) *{\alpha_5},
(28,2) *{\alpha_6},
\ar @{-} "a";"c",
\ar @{-} "c";"d",
\ar @{-} "b";"d",
\ar @{-} "d";"e",
\ar @{-} "e";"f",
\endxy$} & $E_6^6/Sp_4$ & $1,1,1,1,1,1$ & \\
& $E_6({\mathbb C})/E_6$ & $2,2,2,2,2,2$ & \\[1mm]
\hline \rule{0pt}{4mm}
\multirow{2}{*}{\tiny$\xy
\POS (0,0) *\cir<2pt>{} ="a",
(14,5) *\cir<2pt>{} = "b",
(7,0) *\cir<2pt>{}="c",
(14,0) *\cir<2pt>{}="d",
(21,0) *\cir<2pt>{}="e",
(28,0) *\cir<2pt>{}="f",
(35,0) *\cir<2pt>{}="g",
(0,2) *{\alpha_1},
(11,5) *{\alpha_2},
(7,2) *{\alpha_3},
(16,2) *{\alpha_4},
(21,2) *{\alpha_5},
(28,2) *{\alpha_6},
(35,2) *{\alpha_7},
\ar @{-} "a";"c",
\ar @{-} "c";"d",
\ar @{-} "b";"d",
\ar @{-} "d";"e",
\ar @{-} "e";"f",
\ar @{-} "f";"g",
\endxy$}
 & $E_7^7/SU_8$ & $1,1,1,1,1,1,1$ & \\
& $E_7({\mathbb C})/E_7$ & $2,2,2,2,2,2,2$ & \\[1mm]
\hline \rule{0pt}{4mm}
\multirow{2}{*}{\tiny$\xy
\POS (0,0) *\cir<2pt>{} ="a",
(14,5) *\cir<2pt>{} = "b",
(7,0) *\cir<2pt>{}="c",
(14,0) *\cir<2pt>{}="d",
(21,0) *\cir<2pt>{}="e",
(28,0) *\cir<2pt>{}="f",
(35,0) *\cir<2pt>{}="g",
(42,0) *\cir<2pt>{}="h",
(0,2) *{\alpha_1},
(11,5) *{\alpha_2},
(7,2) *{\alpha_3},
(16,2) *{\alpha_4},
(21,2) *{\alpha_5},
(28,2) *{\alpha_6},
(35,2) *{\alpha_7},
(42,2) *{\alpha_8},
\ar @{-} "a";"c",
\ar @{-} "c";"d",
\ar @{-} "b";"d",
\ar @{-} "d";"e",
\ar @{-} "e";"f",
\ar @{-} "f";"g",
\ar @{-} "g";"h"
\endxy $} 
& $E_8^8/SO_{16}$ & $1,1,1,1,1,1,1,1$ & \\
& $E_8({\mathbb C})/E_8$ & $2,2,2,2,2,2,2,2$ & \\[1mm]
\hline \rule{0pt}{4mm}
\multirow{5}{*}{\tiny$\xy
\POS (0,0) *\cir<2pt>{} ="a",
(7,0) *\cir<2pt>{}="b",
(14,0) *\cir<2pt>{}="c",
(21,0) *\cir<2pt>{}="d",
(0,-3) *{\alpha_1},
(7,-3) *{\alpha_2},
(14,-3) *{\alpha_3},
(21,-3) *{\alpha_4},
\ar @{-} "a";"b",
\ar @2{->} "b";"c",
\ar @{-} "c";"d"
\endxy$} & $F_4^4/Sp_3Sp_1$ & $1,1,1,1$ & \\
& $E_6^2/SU_6Sp_1$ & $1,1,2,2$& \\
& $E_7^{-5}/SO_{12}Sp_1$ & $1,1,4,4$ & \\
& $E_8^{-24}/E_7Sp_1$ & $1,1,8,8$ & \\
& $F_4({\mathbb C})/F_4$ & $2,2,2,2$&  \\[1mm]
\hline \rule{0pt}{4mm}
\multirow{2}{*}{\tiny$\xy
\POS (0,0) *\cir<2pt>{} ="a",
(7,0) *\cir<2pt>{}="b",
(0,-3) *{\alpha_1},
(7,-3) *{\alpha_2},
\ar @3{<-} "a";"b"
\endxy$} & $G_2^2/SO_4$ & $1,1$ & \\
& $G_2({\mathbb C})/G_2$ & $2,2$ & \\[1mm]
\hline
\end{tabular}}
\end{table}

Parabolic subalgebras (resp.\ subgroups) of real semisimple Lie algebras (resp.\ Lie groups) play an important role for the geometry of Riemannian symmetric spaces of noncompact type for which their is no analogue in the compact case. We will now describe how to construct all parabolic subalgebras of $\gg$. We denote by $\Psi^+$ the set of positive roots in $\Psi$ with respect to the set $\Lambda$ of simple roots. Let $\Phi$ be a subset of $\Lambda$. We denote by $\Psi_\Phi$ the root subsystem of $\Psi$ generated by $\Phi$, that is, $\Psi_\Phi$ is the intersection of $\Psi$ and the linear span of $\Phi$. We define a reductive subalgebra $\gl_\Phi$ and a nilpotent subalgebra $\gn_\Phi$ of ${\mathfrak g}$  by
\[
\gl_\Phi = \gg_0 \oplus \left(\bigoplus_{\alpha \in \Psi_\Phi} \gg_{\alpha}\right)
\ \ \textrm{and}\ \
\gn_\Phi = \bigoplus_{\alpha \in \Psi^+\setminus \Psi_\Phi^+} {\mathfrak g}_{\alpha}.
\]
It follows from properties of root spaces that $[\gl_\Phi,\gn_\Phi] \subset \gn_\Phi$ and therefore
\[
\gq_\Phi = \gl_\Phi \oplus \gn_\Phi
\]
is a subalgebra of $\gg$, the so-called parabolic subalgebra of $\gg$ associated with the subsystem $\Phi$ of $\Psi$. The decomposition $\gq_\Phi = \gl_\Phi \oplus \gn_\Phi$ is the Chevalley decomposition of the parabolic subalgebra $\gq_\Phi$.

Every parabolic subalgebra of $\gg$ is conjugate in $\gg$ to $\gq_\Phi$ for some subset $\Phi$ of $\Lambda$. The set of conjugacy classes of parabolic subalgebras of $\gg$ therefore has $2^r$ elements. Two parabolic subalgebras $\gq_{\Phi_1}$ and
$\gq_{\Phi_2}$ of $\gg$ are conjugate in the full automorphism group $\Aut(\gg)$ of $\gg$ if and only if there exists an automorphism $F$ of the Dynkin diagram associated to $\Lambda$ with $F(\Phi_1) = \Phi_2$. If $|\Phi| = r-1$ then $\gq_\Phi$ is said to be a maximal parabolic subalgebra of $\gg$.

Let
\[
\ga_\Phi = \bigcap_{\alpha \in \Phi} \ker(\alpha) \subset \ga
\]
be the split component of $\gl_\Phi$ and denote by $\ga^\Phi = \ga \ominus \ga_\Phi$ the orthogonal complement of $\ga_\Phi$ in $\ga$. The reductive subalgebra $\gl_\Phi$ is the centralizer (and the normalizer) of $\ga_\Phi$ in $\gg$.
The orthogonal complement  $\gm_\Phi  = \gl_\Phi \ominus \ga_\Phi$ of $\ga_\Phi$ in $\gl_\Phi$ is a reductive subalgebra of $\gg$. The decomposition
\[
\gq_\Phi = \gm_\Phi \oplus \ga_\Phi \oplus \gn_\Phi
\]
is the Langlands decomposition of the parabolic subalgebra $\gq_\Phi$. We have $[\gm_\Phi,\ga_\Phi] = 0$ and $[\gm_\Phi,\gn_\Phi] \subset \gn_\Phi$. Moreover, $\gg_\Phi = [\gm_\Phi,\gm_\Phi] =[\gl_\Phi,\gl_\Phi]$ is a semisimple subalgebra of $\gg$. The center $\gz_\Phi$ of $\gm_\Phi$ is contained in $\gk_0$ and induces the direct sum decomposition $\gm_\Phi = \gz_\Phi \oplus \gg_\Phi$ and therefore, since $\gz_\Phi \subset \gk_0$, we see that $\gg_\Phi \cap \gk_0 = \gk_0 \ominus\gz_\Phi$.

For each $\alpha \in \Psi$ we define $\gk_\alpha = \gk \cap (\gg_{-\alpha} \oplus \gg_\alpha)$ and $\gp_\alpha = \gp \cap (\gg_{-\alpha} \oplus \gg_\alpha)$. Then we have $\gk_{-\alpha} = \gk_\alpha$, $\gp_{-\alpha} = \gp_\alpha$ and $\gk_\alpha \oplus \gp_\alpha = \gg_{-\alpha} \oplus \gg_\alpha$ for all $\alpha \in \Psi$. From general root space properties it follows that
\[
\gf_\Phi = \gl_\Phi \cap \gp = \ga \oplus \left( \bigoplus_{\alpha \in \Psi_{\Phi}} \gp_\alpha \right)
\ \textrm{and}\
\gb_\Phi = \gm_\Phi \cap \gp = \gg_\Phi \cap \gp = \ga^\Phi \oplus \left( \bigoplus_{\alpha \in \Psi_{\Phi}} \gp_\alpha \right)
\]
are Lie triple systems in $\gp$. We define a subalgebra $\gk_\Phi$ of $\gk$ by
\[
\gk_\Phi = \gq_\Phi \cap \gk = \gl_\Phi \cap \gk = \gm_\Phi \cap\gk = \gk_0 \oplus \left( \bigoplus_{\alpha \in \Psi_{\Phi}}
\gk_\alpha \right).
\]
Then we have $[\gk_\Phi , \gm_\Phi ] \subset \gm_\Phi$, $[\gk_\Phi , \ga_\Phi ] = \{0\}$ and $[\gk_\Phi , \gn_\Phi ] \subset \gn_\Phi$. Moreover, $\gg_\Phi = (\gg_\Phi \cap \gk_\Phi) \oplus \gb_\Phi$ is a Cartan decomposition of the semisimple subalgebra $\gg_\Phi$ of $\gg$ and $\ga^\Phi$ is a maximal abelian subspace of $\gb_\Phi$. If we define $(\gg_\Phi)_0 = (\gg_\Phi \cap \gk_0) \oplus \ga^\Phi$, then $\gg_\Phi = (\gg_\Phi)_0 \oplus \left(\bigoplus_{\alpha \in \Psi_\Phi} \gg_{\alpha}\right)$ is the restricted root space decomposition of $\gg_\Phi$ with respect to ${\mathfrak a}^\Phi$ and $\Phi$ is the corresponding set of simple roots. 

Let $F_\Phi$ and $B_\Phi$ be the connected complete totally geodesic submanifold of $M$ corresponding to the Lie triple systems $\gf_\Phi$ and $\gb_\Phi$, respectively. Then $B_\Phi$ is a Riemannian symmetric space of noncompact type with $\rk(B_\Phi) = |\Phi|$, also known as a boundary component in the maximal Satake compactification of $M$ (see \cite{BJ}). Note that $B_\Phi$ is irreducible if and only if the Dynkin diagram corresponding to $\Phi$ is connected. The totally geodesic submanifold $F_\Phi$ is isometric to the Riemannian product $B_\Phi \times \bbr^{r-|\phi|}$, where $\bbr^{r-|\phi|}$ is the totally geodesic Euclidean space in $M$ corresponding to the abelian Lie triple system $\ga_\Phi$. For $i \in \{1,\ldots,r\}$ we define $\Phi_i = \Lambda \setminus \{\alpha_i\}$, $\gl_i = \gl_{\Phi_i}$, $F_i = F_{\Phi_i}$, $B_i = B_{\Phi_i}$, etcetera. Then we have $F_i = \bbr \times B_i$.

\section{Reflective submanifolds} \label{reflective}

Let $\Sigma'$ be a connected totally geodesic submanifold of $M$. Since $M$ is homogeneous we can assume that $p \in \Sigma'$. Moreover, since every connected totally geodesic submanifold of a Riemannian symmetric space is contained in a connected complete totally geodesic submanifold, we can also assume that $\Sigma'$ is complete. Since $M$ is of noncompact type, $\Sigma'$ is the Riemannian product of a (possibly $0$-dimensional) Euclidean space and a (possibly $0$-dimensional) Riemannian symmetric space of noncompact type. This implies in particular that $\Sigma'$ is simply connected.

The tangent space $\gm' = T_p\Sigma' \subset T_pM = \gp$ is a Lie triple system in $\gp$ and thus $\gg' = [\gm',\gm'] \oplus \gm' \subset \gk \oplus \gp = \gg$ is a subalgebra of $\gg$. Let $G'$ be the connected closed subgroup of $G$ with Lie algebra $\gg'$. Then $\Sigma'$ is the orbit $G' \cdot p$ of the $G'$-action on $M$ containing $p$. Thus we can write $\Sigma' = G'/K'$, where $K' = G'_p$ is the isotropy group of $G'$ at $p$. Since $\Sigma'$ is simply connected, the isotropy group $K'$ is connected. The Lie algebra $\gk'$ of $K'$ is given by $\gk' = [\gm',\gm']$. Note that $G'$ is a  normal subgroup of $G^{\Sigma'} = \{g\in G \mid g(\Sigma') = \Sigma'\}$  and $K'$ is a normal subgroup of $(G^{\Sigma'})_p$. 

The following Slice Lemma was proved in \cite{BO} and will be used later. We formulate it here for the noncompact case, but it is valid also for the compact case. 

\begin {lema} {\sc (Slice Lemma)}\label {Simons} 
Let $M=G/K$ be an irreducible Riemannian symmetric space of noncompact type with $\rk(M)\geq 2$, where $G= I^o(M)$ and $K = G_p$ is the isotropy group of $G$ at $p \in M$. Let $\mathfrak g = \mathfrak k \oplus \mathfrak p$ be the corresponding Cartan decomposition. Let $\Sigma'$ be a nonflat totally geodesic submanifold of $M$ such that  $p\in \Sigma'$.  Let $G'$ be the connected closed subgroup of $G$ with Lie algebra $[\gm',\gm'] \oplus \gm' $, where $T_p\Sigma' = \gm'  \subset \mathfrak p = T_pM$, and $K' = G'_p$. Then the slice representation of $K'$ on $\nu _p \Sigma'$ is nontrivial. 
\end {lema}

In general, the orthogonal complement $\gm''$ of a Lie triple system $\gm'$ in $\gp$ is not a Lie triple system. If $\gm''$ is a Lie triple system, then $\gm'$ is said to be a reflective Lie triple system and $\Sigma'$ is said to be a reflective submanifold of $M$. The notion comes from the fact that the geodesic reflection of $M$ in $\Sigma'$ is a well-defined global isometry of $M$ if and only if both $\gm'$ and $\gm''$ are Lie triple systems. Reflective submanifolds therefore always come in pairs $\Sigma'$ and $\Sigma''$ corresponding to the two reflective Lie triple systems $\gm'$ and $\gm''$. In this situation we write $\Sigma'' = G''/K''$, where $G''$ is the connected closed subgroup of $G$ with Lie algebra $\gg'' = [\gm'',\gm''] \oplus \gm''$ and $K'' = G''_p$ is the connected closed subgroup of $K$ with Lie algebra $\gk'' = [\gm'',\gm'']$. The reflective submanifolds of irreducible simply connected Riemannian symmetric spaces of compact type were classified by Leung (\cite{L1},\cite{L2}). Using duality one obtains the classification of reflective submanifolds in  irreducible Riemannian symmetric spaces of noncompact type.

Let $R$ denote the Riemannian curvature tensor of $M$. As $\Sigma'$ is totally geodesic in $M$, the restriction of $R$ to $\Sigma'$ coincides with the Riemannian curvature tensor of $\Sigma'$. We will regard, via the isotropy representation at $p$, $K' \subset K \subset SO(T_pM)$. Note that $\gk$ and $\gk'$ are generated by the curvature transformations $R_{x,y} \in \so(T_pM)$ with $x,y \in T_pM$ and $x,y \in T_p\Sigma'$, respectively. The curvature operator $\tilde R : \Lambda ^2 (T_pM)  \to \Lambda ^2 (T_pM)$ is negative semidefinite. This implies, as  is well-known,  that $K'$ acts almost effectively on $T_p\Sigma'$.

Let $\rho : K' \to SO (\nu _p\Sigma'),\  k \mapsto d_pk_{\vert \nu _p\Sigma'}$ be the slice representation of $K'$ on the normal space $\nu _p\Sigma'$ of $\Sigma'$ at $p$ and denote by $\ker(\rho)$ the kernel of $\rho$. Let $\chi : K'' \to SO (T _p\Sigma''),\  k \mapsto d_pk_{\vert T _p\Sigma''}$ be the isotropy representation of $K''$ on the tangent space $T_p\Sigma''$.

\begin{lema} \label{ref1}
Let $M = G/K$ be an irreducible Riemannian symmetric space of noncompact type, $\Sigma' = G'/K'$ be a reflective submanifold of $M$ and $\Sigma'' = G''/K''$ be the reflective submanifold of $M$ with $T_p\Sigma'' = \nu_p\Sigma'$. Then:
\begin{itemize}[leftmargin=.3in]
\item[{\rm (i)}] $\rho (K')$ is a normal subgroup of $\chi(K'')$.
\item[{\rm (ii)}] The subspace $(\nu_p\Sigma')_o = \{ \xi \in \nu_p\Sigma' \mid \rho(k')\xi = \xi\ \textrm{for\ all}\ k' \in K'\}$ of $\nu_p\Sigma' = T_p\Sigma''$ is $\chi(K'')$-invariant and $\Sigma'' = \Sigma''_o \times \Sigma''_1$ (Riemannian product), where $\Sigma''_o$ is the totally geodesic submanifold of $\Sigma''$ with $T_p\Sigma''_o = (\nu_p\Sigma')_o$. Moreover, if $\rk(M) \geq 2$, then $\Sigma''_o$ is flat.
\end{itemize}
\end{lema}

\begin{proof}
As previously observed, $K'$ is a normal subgroup of $(G^{\Sigma'})_p$. Observe also that $K'' \subset (G^{\Sigma '})_p$ and that $\rho (K') \subset \chi (K'')$ (for the last inclusion see the paragraph below Lemma 2.1 in \cite{BO}). Then $\rho (K') = \rho (k'' K ' (k'')^{-1}) = \chi (k'') \rho(K') (\chi (k''))^{-1}$ for all $k''\in K''$ and thus $\rho (K')$ is a normal subgroup of $\chi (K'')$. 
Thus the subspace $(\nu_p\Sigma')_o$ of $T_p\Sigma''$ is  $\chi(K'')$-invariant and hence also invariant under the holonomy group of $\Sigma''$ at $p$. Since $\Sigma''$ is simply connected, the de Rham decomposition theorem for Riemannian manifolds implies that $\Sigma''$ decomposes as a Riemannian product into $\Sigma'' = \Sigma''_o \times \Sigma''_1$, where $\Sigma''_o$ is the totally geodesic submanifold of $\Sigma''$ with $T_p\Sigma''_o = (\nu_p\Sigma')_o$.

We write $\Sigma''_o = G''_o/K''_o$, where $G''_o$ is the connected closed subgroup of $G$ with Lie algebra $\gg''_o = [T_p\Sigma''_o,T_p\Sigma''_o] \oplus T_p\Sigma''_o$ and $K''_o$ is the isotropy group of $G''_o$ at $p$. Let $x_1,x_2  \in T_p\Sigma''_o = (\nu_p\Sigma')_o$. For all $y\in T_p\Sigma''_1$ we have $R_{x_1,x_2}y = 0$ since $\Sigma'' = \Sigma''_o \times \Sigma''_1$ is a Riemannian product and totally geodesic in $M$.  Clearly,  $T_p\Sigma''_1$ is $K''_o$-invariant and hence $T_p\Sigma'$ is also $K''_o$-invariant. If $x', y' \in T_p\Sigma'$, then $\langle R_{x_1,x_2}x', y'\rangle = \langle R_{x', y'}x_1,x_2\rangle = 0$, since $x_1, x_2 \in (\nu_p\Sigma')_o$ are fixed under the slice representation of $K'$. Since $\nu _p \Sigma''_o = T_p\Sigma''_1 \oplus T_p\Sigma'$ and $\gk''_o$ is linearly spanned by the curvature endomorphisms of pairs of  elements in $T_p\Sigma''_o$, we conclude that the slice representation of $K''_o$ on $\nu _p \Sigma''_o$ is trivial. It follows from the Slice Lemma \ref{Simons} that $\Sigma''_o$ is flat if $\rk(M) \geq 2$. This finishes the proof of part (ii).
\end{proof}
 
\begin{cor} \label{ref2}
Let $M = G/K$ be an irreducible Riemannian symmetric space of noncompact type with $\rk(M) \geq 2$ and let $\Sigma$ be a totally geodesic submanifold of $M$ which decomposes into a Riemannian product $\Sigma = \Sigma_0 \times \Sigma_1$ with a Euclidean factor $\Sigma_0$ and a semisimple factor $\Sigma_1$ with $\dim(\Sigma_0) > 0$ and $\dim(\Sigma_1) > 0$. Then $\Sigma_1$ is not a reflective submanifold of $M$.
\end{cor}

\begin{proof}
Assume that $\Sigma_1$ is a reflective submanifold of $M$. We will apply Lemma \ref{ref1} with $\Sigma' = \Sigma_1$. In the notation of Lemma \ref{ref1}, we have $T_p\Sigma_0 \subset (\nu_p\Sigma')_o$ and therefore $\Sigma_0$ is contained in the flat factor $\Sigma''_o$ of $\Sigma''$. This implies that $R_{x_0,x''} = 0$ for all $x_0\in T_p\Sigma _0$ and $x'' \in T_p\Sigma''$.
We obviously also have $R_{x_0,x_1} = 0$ for all $x_0 \in T_p\Sigma_0$ and $x_1 \in T_p\Sigma_1 = T_p\Sigma'$. Since $T_pM = T_p\Sigma' \oplus T_p\Sigma''$ this implies $R_{x_0, \cdot } = 0$ for all $x_0 \in T_p\Sigma_0$, which is a contradiction.
\end{proof}

The next result provides a useful sufficient criterion for a semisimple totally geodesic submanifold of an irreducible Riemannian symmetric space to be reflective.

\begin{prop} \label{ref3}
Let $M = G/K$ be an irreducible Riemannian symmetric space of noncompact type with $\rk(M) \geq 2$ and let $\Sigma = G'/K'$ be a semisimple totally geodesic submanifold of $M$. Assume that the kernel $\ker(\rho)$ of the  slice representation $\rho : K' \to SO (\nu _p \Sigma)$ has positive dimension. Then we have
\[
\nu _p \Sigma = \{\xi \in T_pM \mid \rho(k)\xi = \xi\ \textrm{for\ all}\ k \in \ker(\rho)^o\}
\]
and, in particular, $\Sigma$ is a reflective submanifold of $M$.
\end{prop}

\begin{proof} 
The subspace $\bbv =  \{\xi \in T_p\Sigma \mid d_pk(\xi) = \xi\ \textrm{for\ all}\ k \in \ker(\rho)^o\}$ of $T_p\Sigma$ is $K'$-invariant since $\ker(\rho)^o$ is a normal subgroup of $K'$.

We first assume  that $\mathbb V = T_p \Sigma$. Since $\ker(\rho)^o$ acts trivially on $\nu _p\Sigma$ we conclude that $\ker(\rho)^o$ and hence $K'$ acts trivially on $T_pM$, which is a contradiction.

Next, we  assume that  $\bbv$ is a nontrivial proper $K'$-invariant  subspace of $T_p\Sigma$. Then $\Sigma$ decomposes as a Riemannian product $\Sigma = \Sigma_1 \times \Sigma_2$, where $\bbv = T_p\Sigma_1$. If we write, as usual, $\Sigma_i = G_i/K_i$, then $K' = K_1 \times K_2$ (almost direct product). Let $\gh _i$  be the orthogonal projection of the Lie algebra of $\ker(\rho)$ into $\gk_i$ and let $H_i$ be the connected subgroup of $K_i$ with Lie algebra $\gh_i$. Then $H_1$ acts trivially on $\bbv = T_p\Sigma_1$ since both $\ker(\rho)^o$ and $H_2$ act trivially on $\bbv$. Since $K_1$ acts almost effectively on $T_p\Sigma_1$ and $H_1$ is connected, it follows that $H_1$ is trivial. Thus we have shown that $\ker(\rho)^o \subset K_2$.

Note that  $\{\xi \in T_pM \mid \rho(k)\xi = \xi\ \textrm{for\ all}\ k \in \ker(\rho)^o\} = \mathbb V \oplus \nu _p \Sigma = \nu _p \Sigma_2$.  This shows that $\Sigma _2$ is a reflective submanifold of $M$.  Let $\Sigma _3 = G_3/K_3$ be the reflective submanifold of $M$ with $T_p\Sigma_3 = \nu _p\Sigma _2$. We denote by $\rho _i : K_i \to SO (\nu _p \Sigma _i)$ the slice representation of $K_i$ on the normal space $\nu_p\Sigma_i$ and by $\chi_i : K_i \to SO(T_p\Sigma_i)$ the isotropy representation of $K_i$, $i \in \{1,2,3\}$.

From Lemma \ref{ref1}(i) we see that $\rho_2 (K_2)$ is a normal subgroup of $\chi_3(K_3)$. Let $\mathbb W$ be the set of fixed vectors of $\rho _2 (K_2)$ in $\nu _p \Sigma _2 = T_p \Sigma _3 = \mathbb V \oplus \nu _p \Sigma  = T_p\Sigma _1 \oplus \nu _p \Sigma$. Since $K_2$ acts trivially on $T_p\Sigma _1$ one has that $T_p\Sigma _1 \subset \mathbb W$. From Lemma \ref{ref1}(ii) we know that $\mathbb W$ is the tangent space of a Euclidean factor of $\Sigma _3$. This is a contradiction since 
$\Sigma _1$ is  contained in this Euclidean factor, however,  $\Sigma _1$ is not flat as $\Sigma $ is semisimple. It follows that $\mathbb V =\{ 0 \}$, which proves the assertion.
\end{proof}

The following consequence of Proposition \ref{ref3} states that totally geodesic submanifolds of sufficiently small codimension in irreducible Riemannian symmetric spaces are reflective. 

\begin{cor} \label{ref4} 
Let $M$ be an $n$-dimensional irreducible Riemannian symmetric space of noncompact type with $r = \rk(M) \geq 2$ and let $\Sigma$ be a semisimple connected complete totally geodesic submanifold of $M$ with $\codim(\Sigma) = d$. If \[\frac{1}{2}d(d+1) +\rk(\Sigma) < n,\] then $\Sigma$ is a reflective submanifold of $M$. In particular, if \[d(d+1) < 2(n-r),\] then $\Sigma$ is a reflective submanifold of $M$.
\end{cor}

\begin{proof}
As usual, we write $\Sigma = G'/K'$. If $\dim(K') > \dim(SO(\nu_p\Sigma)) = \frac{1}{2}d(d-1)$, then the kernel of the slice representation $\rho : K' \to SO(\nu_p\Sigma)$ must have positive dimension and therefore $\Sigma$ is a reflective submanifold of $M$ by Proposition \ref{ref3}. A principal $K'$-orbit on $\Sigma$ has dimension $n-d-\rk(\Sigma)$ and thus $\dim(K') \geq n - d - \rk(\Sigma)$. Consequently, if $\frac{1}{2}d(d-1) < n - d - \rk(\Sigma)$, then $\Sigma$ is a reflective submanifold of $M$. The inequality $\frac{1}{2}d(d-1) < n - d - \rk(\Sigma)$ is equivalent to $\frac{1}{2}d(d+1) + \rk(\Sigma) < n$. The last statement follows from the fact that $\rk(\Sigma) \leq \rk(M) = r$.
\end{proof}

\section{Non-semisimple maximal totally geodesic submanifolds} \label{nonsemisimple}

Let $\Sigma$ be a connected totally geodesic submanifold of $M$. We may assume that $\Sigma$ is complete and $p \in \Sigma$. Every connected complete totally geodesic submanifold of a Riemannian symmetric space is again a Riemannian symmetric space. In this paper, when we consider a totally geodesic submanifold $\Sigma$ of $M$, we always assume that $p \in \Sigma$ and that $\Sigma$ is connected and complete. Since $M$ is of noncompact type, it follows that $\Sigma$ is isometric to the Riemannian product $\Sigma_0 \times \Sigma_1$, where $\Sigma_0$ is a (possibly $0$-dimensional) Euclidean space and $\Sigma_1$ is a (possibly $0$-dimensional) Riemannian symmetric space of noncompact type. 

The next result relates non-semisimple maximal totally geodesic submanifolds of $M$ to the reductive factors in the Chevalley decompositions of the maximal parabolic subalgebras of $\gg$.

\begin{prop}\label{totgeodnss}
Let $M = G/K$ be an irreducible Riemannian symmetric space of noncompact type and let $\Sigma$ be a  non-semisimple maximal totally geodesic submanifold of $M$. Then $\Sigma$ is congruent to $F_i = \bbr \times B_i$ for some $i \in \{i,\ldots,r\}$.
\end{prop}

\begin{proof}
Let $\ga$ be a maximal abelian subspace of $\gp$ with $T_p\Sigma_0 \subset \ga$ and consider the restricted root space decomposition of $\gg$ induced by $\ga$. We define $\Upsilon = \{ \alpha_i \in \Lambda \mid \alpha_i(T_p\Sigma_0) = 0\} \subset \Lambda$. Assume that $\Upsilon = \Lambda$, which means that $\alpha_i(T_p\Sigma_0) = 0$ for all $\alpha_i \in \Lambda$. This implies $T_p\Sigma_0 = \{0\}$ and therefore $\Sigma = \Sigma_1$ is semisimple, which is a contradiction. Thus we have $|\Upsilon| < |\Lambda| = r$ and therefore there exists $i \in \{1,\ldots,r\}$ such that $\Upsilon \subset \Phi_i$. Then we get
\begin{eqnarray*}
T_p\Sigma 
& \subset & Z_\gp(T_p\Sigma_0) = \{ X \in \gp \mid [X,Y] = 0\ \textrm{for\ all}\ Y \in T_p\Sigma_0\} \\
& \subset & Z_\gg(T_p\Sigma_0) = \{ X \in \gg \mid [X,Y] = 0\ \textrm{for\ all}\ Y \in T_p\Sigma_0\} \\
& = & \gg_0 \oplus \left( \bigoplus_{\alpha \in \Psi,\alpha(T_p\Sigma_0) = \{0\}} \gg_\alpha \right) 
= \gl_\Upsilon  \subset  \gl_i,
\end{eqnarray*}
which implies $T_p\Sigma \subset \gl_i \cap \gp = \gf_i$ and therefore $\Sigma \subset F_i = \bbr \times B_i$. 
If $\Sigma$ is a maximal totally geodesic submanifold of $M$ we must have $\Sigma = F_i$, since $F_i$ is a totally geodesic submanifold of $M$ which is strictly contained in $M$. 
\end{proof}

The remaining problem is to clarify which of the totally geodesic submanifolds $F_i$ are maximal. The solution of this problem is related to symmetric R-spaces. Let $M = G/K$ be an irreducible Riemannian symmetric space of noncompact type and consider the isotropy representation
\[ \chi : K \to T_pM = \gp,\ v \mapsto d_pk(v) = \Ad(k)v.\]
For every $0 \neq v \in \gp$ the orbit
\[ K\cdot v = \{\Ad(k)v \mid k \in K\} \subset \gp \]
is called an R-space (or real flag manifold). One can show that the normal space $\nu_v(K \cdot v)$ of $K \cdot v$ at $v$ is equal to
\[
\nu_v(K \cdot v) = Z_\gp(v) = \{ w \in \gp \mid [v,w] = 0 \},
\]
where $Z_\gp(v)$ is the centralizer of $v$ in $\gp$. It follows from the Jacobi identity that $Z_\gp(v)$ is a Lie triple system. Thus, for every $0 \neq v \in \gp$, there exists a connected complete totally geodesic submanifold $\Sigma^v$ of $M$ with $T_p\Sigma^v = \nu_v(K \cdot v)$. Since every $v \in \gp$ is contained in a maximal abelian subspace of $\gp$ we can assume that $v \in \ga$. Then we have $\gl_\Phi = Z_\gg(v)$ with $\Phi = \{\alpha_i \in \Lambda \mid \alpha_i(v) = 0\}$, which implies $\gf_\Phi = Z_\gp(v) = \nu_v(K \cdot v)$ and therefore $F_\Phi = \Sigma^v$. 

\begin{table}[b]
\caption{Highest roots $\delta$ of root systems $(R)$} 
\label{highestroot}
{\footnotesize\begin{tabular}{ | p{1cm}  p{6.5cm}  p{1.5cm} |}
\hline \rule{0pt}{4mm}
\hspace{-2mm}  $(R)$ & Highest root $\delta = \delta_1\alpha_1 + \ldots + \delta_r\alpha_r$ & Comments\\[1mm]
\hline \rule{0pt}{4mm}
\hspace{-2mm} 
$(A_r)$ & $\alpha_1 + \ldots + \alpha_r$ & $r \geq 1$\\
$(B_r)$ & $\alpha_1 + 2\alpha_2 + \ldots + 2\alpha_r$ & $r \geq 2$\\
$(C_r)$ & $2\alpha_1 + \ldots + 2\alpha_{r-1} + \alpha_r$ & $r \geq 3$\\
$(D_r)$ & $\alpha_1 + 2\alpha_2 + \ldots + 2\alpha_{r-2} + \alpha_{r-1} + \alpha_r$ & $r \geq 4$\\
$(BC_r)$ & $2\alpha_1 + \ldots + 2\alpha_r$ & $r \geq 1$\\
$(E_6)$ & $\alpha_1 + 2\alpha_2 + 2\alpha_3 + 3\alpha_4 + 2\alpha_5 + \alpha_6$ & \\
$(E_7)$ & $2\alpha_1 + 2\alpha_2 + 3\alpha_3 + 4\alpha_4 + 3\alpha_5 + 2\alpha_6 + \alpha_7$ & \\
$(E_8)$ & $2\alpha_1 + 3\alpha_2 + 4\alpha_3 + 6\alpha_4 + 5\alpha_5 + 4\alpha_6 + 3\alpha_7 + 2\alpha_8$ & \\
$(F_4)$ & $2\alpha_1 + 3\alpha_2 + 4\alpha_3 + 2\alpha_4$ & \\
$(G_2)$ & $3\alpha_1 + 2\alpha_2$ & \\[1mm]
\hline 
\end{tabular}}
\end{table}

A special situation occurs when the orbit $K \cdot v$ is a symmetric space. In this situation the orbit $K \cdot v$ is called an irreducible symmetric R-space. The irreducibility here refers to the irreducibility of the symmetric space $G/K$ and not to the orbit. An irreducible symmetric R-space can be reducible as a Riemannian manifold. The irreducible symmetric R-spaces were classified by Kobayashi and Nagano in \cite{KN}. Their classification can be read off from the Dynkin diagram and highest root of the symmetric spaces $G/K$. In Table \ref{dynkin} we already listed the Dynkin diagrams. In Table \ref{highestroot} we list the corresponding highest roots $\delta = \delta_1\alpha_1 + \ldots + \delta_r\alpha_r$. 

Kobayashi and Nagano proved that an R-space $K\cdot v$ is symmetric if and only if $v = H^i$ and $\delta_i = 1$. From Tables \ref{dynkin} and \ref{highestroot} one can easily get the classification of irreducible symmetric R-spaces. We can now state the main result of this section:

\begin{teo} \label{rspace}
Let $M = G/K$ be an irreducible Riemannian symmetric space of noncompact type and let $\Sigma$ be a non-semisimple connected complete totally geodesic submanifold of $M$. Then the following statements are equivalent:
\begin{itemize}[leftmargin=.3in]
\item[\rm{(i)}] $\Sigma$ is a maximal totally geodesic submanifold of $M$;
\item[\rm{(ii)}] $\Sigma$ is isometrically congruent to $F_i = \bbr \times B_i$ and $\delta_i = 1$;
\item[\rm{(iii)}] $\nu_p\Sigma$ is the tangent space of a symmetric R-space in $T_pM$;
\item[\rm{(iv)}] The pair $(M,\Sigma)$ is as in Table {\rm \ref{totgeodnsstable}}.
\end{itemize}
\end{teo}

\begin{table}[h]
\caption{Non-semisimple maximal totally geodesic submanifolds $\Sigma = \bbr \times B$ of irreducible Riemannian symmetric spaces $M$ of noncompact type} 
\label{totgeodnsstable}
{\footnotesize\begin{tabular}{ | p{2.7cm}  p{4.9cm}  p{1.4cm}  p{2.6cm} | }
\hline \rule{0pt}{4mm}
\hspace{-2mm}  
$M$ & $B$ & $\codim(\Sigma)$ & Comments\\[1mm]
\hline \rule{0pt}{4mm}
\hspace{-2mm}  
$SL_{r+1}(\bbr)/SO_{r+1}$ & $SL_i(\bbr)/SO_i \times SL_{r+1-i}(\bbr)/SO_{r+1-i}$ & $i(r+1-i)$ & $r \geq 2$, $1 \leq i \leq [r/2]$\\
$SL_{r+1}(\bbc)/SU_{r+1}$ & $SL_i(\bbc)/SU_i \times SU_{r+1-i}(\bbc)/SU_{r+1-i}$ & $2i(r+1-i)$ & $r \geq 2$, $1 \leq i \leq [r/2]$\\
$SU^*_{2r+2}/Sp_{r+1}$ & $SU^*_{2i}/Sp_i \times SU^*_{2(r+1-i)}/Sp_{r+1-i}$ & $4i(r+1-i)$ & $r \geq 2$, $1 \leq i \leq [r/2]$\\
$E_6^{-26}/F_4$ & $\bbr H^9$ & $16$ & \\[1mm]
\hline \rule{0pt}{4mm}
\hspace{-2mm}  
$SO^o_{r,r+k}/SO_{r}SO_{r+k}$ & $SO^o_{r-1,r-1+k}/SO_{r-1}SO_{r-1+k}$ & $2r-2+k$ &  $r \geq 2, k \geq 1$\\
$SO_{2r+1}(\bbc)/SO_{2r+1}$ & $SO_{2r-1}(\bbc)/SO_{2r-1}$& $4r-2$ & $r \geq 2$\\[1mm]
\hline \rule{0pt}{4mm}
\hspace{-2mm}  
$Sp_r(\bbr)/U_r$ &  $SL_r(\bbr)/SO_r$ & $\frac{1}{2}r(r+1)$ & $r \geq 3$\\
$SU_{r,r}/S(U_rU_r)$ &  $SL_r(\bbc)/SU_r $ & $r^2$ & $r \geq 3$\\
$Sp_r(\bbc)/Sp_r$ & $SL_r(\bbc)/SU_r$ & $r(r+1)$ &  $r \geq 3$\\
$SO^*_{4r}/U_{2r}$ & $SU^*_{2r}/Sp_r$ & $r(2r-1)$ &  $r \geq 3$\\
$Sp_{r,r}/Sp_rSp_r$ &  $SU^*_{2r}/Sp_r$ & $r(2r+1)$ & $r \geq 2$\\
$E_7^{-25}/E_6U_1$ & $E_6^{-26}/F_4$ & $27$ & \\[1mm]
\hline \rule{0pt}{4mm}
\hspace{-2mm}  
$SO^o_{r,r}/SO_{r}SO_{r}$ & $SO^o_{r-1,r-1}/SO_{r-1}SO_{r-1}$ & $2(r-1)$ & $r \geq 4$\\
& $SL_r(\bbr)/SO_r$ & $\frac{1}{2}r(r-1)$ & $r \geq 4$\\
$SO_{2r}(\bbc)/SO_{2r}$ & $SO_{2(r-1)}(\bbc)/SO_{2(r-1)}$& $4(r-1)$  & $r \geq 4$\\
& $SL_r(\bbc)/SU_r$ & $r(r-1)$ & $r \geq 4$\\[1mm]
\hline \rule{0pt}{4mm}
\hspace{-2mm}  
$E_6^6/Sp_4$ & $SO^o_{5,5}/SO_5SO_5$ &  $16$ & \\
$E_7^7/SU_8$ & $E_6^6/Sp_4$ & $27$ & \\[1mm]
\hline \rule{0pt}{4mm}
\hspace{-2mm}  
$E_6({\mathbb C})/E_6$ & $SO_{10}(\bbc)/SO_{10}$ & $32$ & \\
$E_7({\mathbb C})/E_7$ & $E_6(\bbc)/E_6$ & $54$ & \\[1mm]
\hline
\end{tabular}}
\end{table}

\begin{proof}
The equivalence of (ii) and (iv) is a straightforward computation using Tables \ref{dynkin} and \ref{highestroot}. Kobayashi and Nagano proved that an R-space $K \cdot v$ is symmetric if and only if $v = H^i$ and $\delta_i = 1$. In this situation we have  $\nu_{H^i}(K \cdot H^i) = Z_\gp(H^i) = \gf_i = T_pF_i$ and hence $\nu_pF_i = T_{H^i}(K \cdot H^i)$. This gives the equivalence of (ii) and (iii). We shall now prove the equivalence of (i) and (ii).

We first assume that $\Sigma$ is a maximal totally geodesic submanifold of $M$. From Proposition \ref{totgeodnss} we know that, up to conjugacy, $\Sigma = F_i$ for some $i \in \{1,\ldots,r\}$. Assume that $\delta_i \geq 2$ and let $t$ be a prime number with $t \leq \delta_i$. Then define the semisimple subalgebra $\gh_{i,t}$ of $\gg$ by
\[
\gh_{i,t} = \gg_0 \oplus \left( \bigoplus_{\alpha \in \Psi,\alpha(H^i) \equiv 0(\textrm{mod}\,t)} \gg_\alpha \right) .
\]
Since 
\[
\gl_i = \gg_0 \oplus \left( \bigoplus_{\alpha \in \Psi,\alpha(H^i) = 0} \gg_\alpha \right) 
\]
and $\delta_i \geq t$ we see that $\gl_i$ is strictly contained in $\gh_{i,t}$. It follows that the Lie triple system $\gl_i \cap \gp = \gf_i$ is strictly contained in the Lie triple system $\gh_{i,t} \cap \gp$. This is a contradiction since, by assumption, $T_p\Sigma = \gf_i$ is a maximal Lie triple system. Consequently we must have $\delta_i = 1$. This finishes the proof for ``(i) $\Rightarrow$ (ii)''.

Conversely, let us assume that $\Sigma = F_i$ for some $i \in \{1,\ldots,r\}$ and that $\delta_i = 1$. We denote by $S_i$ the symmetric R-space $K \cdot H^i \subset \gp = T_pM$. Then we have $\alpha(H^i) \in \{-1,0,+1\}$ for all $\alpha \in \Psi$ and therefore $\ad(H^i)^2$ induces a vector space decomposition $\gg = \gg^0 \oplus \gg^1$ of $\gg$, where
\[
\gg^0 = \gl_i = \gg_0 \oplus \left( \bigoplus_{\alpha \in \Psi,\alpha(H^i) = 0} \gg_\alpha \right)\ \textrm{and}\ 
\gg^1 = \bigoplus_{\alpha \in \Psi,\alpha(H^i) = \pm 1} \gg_\alpha .
\]
The map $X_0 + X_1 \to X_0 - X_1$ defines an involutive automorphism of $\gg = \gg^0 \oplus \gg^1$. We denote by $s_i : \gp \to \gp$ the induced isomorphism on $\gp$. Then we have $s_i(X) = -X$ for all $X \in T_{H^i}S_i = \gg^1 \cap \gp = \oplus_{\alpha \in \Psi,\alpha(H^i) = 1} \gp_\alpha $ and $s_i(X) = X$ for all $X \in \nu_{H^i}S_i = \gg^0 \cap \gp = \gf_i = \bbr \times \gb_i = Z_\gp(H^i)$. The isomorphism $s_i$ is the orthogonal reflection of $\gp$ in the normal space $\nu_{H^i}S_i$ and its restriction to $S_i$ leaves $S_i$ invariant and hence induces an involutive isometry on $S_i$ for which $H^i$ is an isolated fixed point. This shows that $S_i$ is a symmetric R-space and that $S_i$ is an extrinsically symmetric submanifold of the Euclidean space $T_pM = \gp$ with $s_i$ as the extrinsic symmetry.

Since $[\gg^\nu,\gg^\mu] \subset \gg^{(\nu+\mu)({\rm mod}\,2)}$, we see that $\nu_{H^i}S_i =  \gg^0 \cap \gp$ and $T_{H^i}S_i = \gg^1 \cap \gp$ are Lie triple systems. It follows that both the tangent space and the normal space of the symmetric R-space $S_i$ at $H^i$ are reflective Lie triple systems.

Let $\bbv \neq \gp$ be a Lie triple system in $\gp$ with $\gf_i \subset \bbv$ and let $\Sigma'$ be the connected complete totally geodesic submanifold of $M$ with $T_p\Sigma' = \bbv$. Then we have $\Sigma' = G'/K'$, where $G'$ and $K'$ is the connected closed subgroup of $G$ with Lie algebra $\gg' = [\bbv,\bbv] \oplus \bbv$ and $\gk' = [\bbv,\bbv]$, respectively. 

Since $\rk(M) \geq 2$, the semisimple factor $\gb_i$ of $\gf_i$ is non-trivial and therefore $\bbv$ is a non-abelian subspace of $\gp$. Since $T_{H^i}S_i$ is a Lie triple system, $\bbv \cap T_{H^i}S_i$ is a Lie triple system as well. Let $N$ be the connected component containing $H^i$ of the intersection $\Sigma' \cap S_i$. It is clear from the construction that $N$ is a smooth submanifold of $S_i$ in an open neighborhood of $H^i$.

We identify $X\in \gg$ with the Killing field $q\mapsto X.q = \frac{d}{dt}_{|t = 0}(t \mapsto \Exp(tX)(q))$ on $M$. The orthogonal projection $\bar X$  of $X_{\vert \Sigma'}$ to $T\Sigma'$ is a Killing field on the totally geodesic submanifold $\Sigma'$ which lies in the transvection algebra of $\Sigma'$ (see the paragraph below Lemma 2.1 in \cite {BO}). Note that $\bar X.p=0$ if $X\in \gk$. Then, if $X\in \gk$, there exists $X' \in \gk'$ such that $Z_{\vert \Sigma'}$ is always perpendicular to $\Sigma'$, where $Z=X-X'$. This implies that $Z.\bbv\subset \bbv^\perp$.  In fact, let $\gamma _u$ be the geodesic in $M$ with initial condition $\gamma_u'(0) = u \in \bbv$. The Jacobi field $Z.\gamma _u(t)$ is perpendicular to $T_{\gamma_u(t)} \Sigma'$ and therefore its covariant derivative $(Z.\gamma _u)' (t)$ must be so. Hence $(Z.\gamma _u)'(0)= Z.u \in \bbv ^\perp$. So, if $X\in \gk$, we have $X.u = X'.u  + Z.u$ and thus $T_u(K\cdot u) \subset   T_u(K'\cdot u) \oplus \bbv^\perp$ for all $u \in \bbv$. This implies $T_u(K\cdot u) = T_u(K'\cdot u)\oplus \bbv^\perp$ (orthogonal direct sum) and  $\nu _u (K'\cdot u) = \nu _u (K\cdot u)$ for all $u\in S_i$, since $\nu _{H^i}S_i \subset \bbv$.

As we have previously observed, $N$ is a submanifold of $S_i$ in an open neighborhood of $H^i$. Since $K'\cdot H^i \subset \bbv$ and $K'\subset K$ we obtain $K'\cdot H^i \subset N$ and $K'\cdot N = N$. From the previous paragraph we conclude that $N$ coincides with $K'\cdot H^i$ around $H^i$, since both submanifolds of $S_i$ have the same dimension. Moreover, since $\bbv$ is $K'$-invariant, $\bbv$ contains the normal space $\nu _w S_i= \nu _w (K\cdot w)$ for all $w \in  K'\cdot H^i$. This implies in particular that $N$ is totally geodesic in $S_i$ at all points $w \in  K'\cdot H^i$. Thus $N$ is a submanifold around any $w\in K'\cdot H^i$ and $N$ coincides, around $w$, with this orbit. Therefore $K'\cdot H^i$ is an open subset of $N$. Since $K'\cdot H^i$ is compact and $N$ is Hausdorff, the orbit $K'\cdot H^i$ is a closed subset of $N$. Since $N$ is connected this implies that $N = K'\cdot H^i$ is a totally geodesic submanifold of $S_i$. 

Let us consider the extrinsic symmetry $s_i$ at $H^i$ of the extrinsically symmetric submanifold $S_i$ of $T_pM$. Since $s_i$ leaves  $S_i$,  $\bbv$ and $\{H^i\}$ invariant, it also leaves $N = K'\cdot H^i$, the connected component of $S_i \cap \bbv$ containing $H^i$, invariant. Hence $s_i$ restricted to $\bbv$ is an extrinsic symmetry of $N\subset \bbv$ at $H^i$. This proves that $N$ is an extrinsically symmetric submanifold of $\mathbb V$.

Note that the extrinsic symmetry $s_i$ has the property $s_i (\bbv) = \bbv$ and therefore $s_i K's_i^{-1} = K'$. 

We want to prove that $N = \{ H^i\}$, or equivalently, that $\bbv = \nu_{H^i}S_i$. Assume that this is not true. Let $\bbw \subset \bbv$ be the linear span of $N = K'\cdot H^i$. Then $\bbw$ is the tangent space to a (symmetric) Riemannian factor of $\Sigma'$, since it is $K'$-invariant. The subspace $\bbw$ cannot have an abelian part since $N = K'\cdot H^i$ is full in $\bbw$.  Also, since $K$ acts irreducibly on $T_pM$, $K$ must act effectively on the symmetric orbit $S_i$. The group $K$ is generated by the so-called geometric transvections $\{ s_x\circ s_y\} $, where $x, y \in S_i$ and $s_x$ denotes the extrinsic symmetry at $x$. In fact, the connected group $K$ cannot be bigger than the group of transvections of the symmetric space $S_i$ since $S_i$ is compact, and so any Killing field on $S_i$ is bounded and hence belongs to the Lie algebra of the transvection group. 

Let $K''$ be the connected closed subgroup of $K'$ with Lie algebra $\gk'' =  [\bbw,\bbw] \oplus \bbw$. Note that $K'\cdot H^i = K''\cdot H^i$. Moreover, $K''$ acts almost effectively on $N$. In fact, $K''$ acts almost effectively on $\bbw$ (see Section 2 of \cite {BO}) and if $k''\in K''$ acts trivially on $N$ it must act trivially on its linear span. We also have  $K ' = K'' \times \bar{K}$ (almost direct product), where $\bar{K}$ is the connected closed subgroup of $K'$ with Lie algebra $\bar{\gk} =  [\bbw^\perp \cap \bbv,\bbw^\perp \cap \bbv] \oplus (\bbw^\perp \cap \bbv)$. 

As we have seen above,  $N = K''\cdot H^i$ is a symmetric submanifold of $S_i$ and thus $K''$ must be  generated by $\{ s_{x '}\circ s_{y'}\}$ with $x', y' \in K'\cdot H^i$. The following observation is crucial: $\{ s_{x '}\circ s_{y'}\} $ is the identity on the orthogonal complement of $\bbw$. In fact, $s_{x '}$ is the identity on $\bbw^\perp \cap \bbv$, since this subspace is contained in $\nu _{x'} S_i$. Moreover, $s_{x '}$ is minus the identity on $\bbv^\perp$, which is tangent to $S_i$ at $x'$. The same is true if one replaces $x'$ by $y'$ and so $s_{x '}\circ s_{y'}$ is the identity on $\bbv^\perp \oplus (\bbw^\perp \cap \bbv) = \bbw^\perp $. This implies that $K''$ acts trivially on $\bbw^\perp$, which contradicts the Slice Lemma \ref{Simons}. Then $\bbv = \nu_{H^i}S_i$ which implies that $\nu_{H^i} S_i = T_p\Sigma$ is maximal. Thus we have proved that $\Sigma$ is a maximal totally geodesic submanifold of $M$. This finishes the proof of ``(ii) $\Rightarrow$ (i)''.
\end{proof}

From Theorem \ref{rspace} and Table \ref{highestroot} we obtain

\begin{cor} \label{nonss}
Let $M = G/K$ be an irreducible Riemannian symmetric space of noncompact type. If the restricted root system of $M$ is of type $(BC_r)$, $(E_8)$, $(F_4)$ or $(G_2)$, then every maximal totally geodesic submanifold of $M$ is semisimple.
\end{cor}

We have seen in the proof of Theorem \ref{rspace} that $\nu_pF_i$ is a Lie triple system when $\delta_i = 1$, which implies that $T_pF_i$ is a reflective Lie triple system when $\delta_i = 1$. From Theorem \ref{rspace} we can therefore conclude:

\begin{cor} \label{nssrefl}
Let $M = G/K$ be an irreducible Riemannian symmetric space of noncompact type and let $\Sigma$ be a non-semisimple maximal totally geodesic submanifold of $M$. Then $\Sigma$ is a reflective submanifold of $M$.
\end{cor}

We remark that the analogous statement for the semisimple case does not hold. For example, $SL_3(\bbr)/SO_3$ is a semisimple maximal totally geodesic submanifold of $G_2^2/SO_4$ which is not reflective. 

We recall from \cite{BO} the following result:

\begin{teo}\label{bound}
Let $M$ be an irreducible Riemannian symmetric space. Then 
\[\rk(M) \leq i(M).\]
\end{teo}

From Table \ref{totgeodnsstable} we obtain that the codimension of the totally geodesic submanifold $\Sigma = \bbr \times  SL_{r}(\bbr)/SO_{r}$ in $M =  SL_{r+1}(\bbr)/SO_{r+1}$ is equal to $r = \rk(M)$, which implies $i(M) \leq \rk(M)$. Using Theorem \ref{bound} we thus conclude:

\begin{cor} \label{SLSO}
For $M = SL_{r+1}(\bbr)/SO_{r+1}$ we have $\rk(M) = r = i(M)$.
\end{cor}

\section{Examples of symmetric spaces with $\rk(M) = i(M)$} \label{ex}

We first consider the symmetric space $M = SL_{r+1}({\mathbb R})/SO_{r+1}$ for $r \geq 1$ and present a more explicit version of Corollary \ref{SLSO}.  This symmetric space has $\rk(M) = r$ and $\dim(M) = \frac{1}{2}r(r+3)$. For $r = 1$ we get the real hyperbolic plane $\bbr H^2$. Thus, if $\Sigma$ is a geodesic in $M$, we have $\codim(\Sigma) = 1 = \rk(M)$. For $r \geq 2$ we consider the Cartan decomposition $\gg = \gk \oplus \gp$ of the Lie algebra $\gg= \gsl_{r+1}({\mathbb R})$ of $G = SL_{r+1}({\mathbb R})$ which is induced by the Lie algebra ${\mathfrak k} = \so_{r+1}$ of $K = SO_{r+1}$. The vector space $\gp$ is given by 
\[ 
\gp = \{A \in \gsl_{r+1}({\mathbb R}) \mid A^T = A\}.
\]
We now define
\[ 
\gm = \left\{ \left. A = \begin{pmatrix} -\tr(B) & 0 \\ 0 & B \end{pmatrix} \in \gp\  \right| \ B \in \ggl_r(\bbr),\ B^T = B \right\} .
\]
Then we have
\[ 
[[\gm,\gm],\gm] = \left\{ \left. A = \begin{pmatrix} 0 & 0 \\ 0 & \ B \end{pmatrix} \in \gp \ \right| \ B \in \gsl_r(\bbr),\ B^T = B \right\} \subset \gm,
\]
which shows that $\gm$ is a Lie triple system in $\gp$. We have $\dim(\gm) = \frac{1}{2}r(r+1)$ and hence $\dim(\gp) - \dim(\gm) = \frac{1}{2}r(r+3) - \frac{1}{2}r(r+1) = r$. Thus the connected complete totally geodesic submanifold $\Sigma$ of $M$ corresponding to the Lie triple system ${\mathfrak m}$, which is isometric to $\bbr \times SL_r({\mathbb R})/SO_r$, satisfies $\codim(\Sigma) = r = \rk(M)$. From Theorem \ref{bound} we can therefore conclude that the index of $SL_{r+1}({\mathbb R})/SO_{r+1}$ is equal to the rank of $SL_{r+1}({\mathbb R})/SO_{r+1}$. We remark that $\bbr \times SL_r({\mathbb R})/SO_r$ is tangent to the normal space of a Veronese embedding of the real projective space $\bbr P^r$ into $\gp$ (see e.g.\ Lemma 8.1 in \cite{ORi}).

Next, we consider the symmetric space $M = SO^o_{r,r+k}/SO_rSO_{r+k}$ with $r \geq 1$, $k \geq 0$ and $(r,k) \notin \{(1,0),(2,0)\}$. This symmetric space has $\rk(M) = r$ and $\dim(M) = r(r+k)$. For $(r,k) = (1,0)$ we have $\dim(M) = 1$ and so $M$ is not of noncompact type. For $(r,k) = (2,0)$ we have the symmetric space $M = SO^o_{2,2}/SO_2SO_2$ which is isometric to the Riemannian product of two real hyperbolic planes and therefore not irreducible. Note that $SO^o_{1,2}/SO_2 = SL_2({\mathbb R})/SO_2$ and $SO^o_{3,3}/SO_3SO_3 = SL_4({\mathbb R})/SO_4$.

For $r = 1$ we get the $(k+1)$-dimensional real hyperbolic space $M = \bbr H^{k+1} = SO^o_{1,1+k}/SO_{1+k}$. This space contains a totally geodesic hypersurface $\Sigma = \bbr H^k$ and therefore $\rk(M) = 1 = i(M)$. 

Now assume that $r \geq 2$ and consider the Cartan decomposition $\gg = \gk \oplus \gp$ of the Lie algebra $\gg= \gso_{r,r+k}$ of $G = SO^o_{r,r+k}$ which is induced by the Lie algebra ${\mathfrak k} = \gso_r \oplus \so_{r+k}$ of $K = SO_rSO_{r+k}$. The vector space $\gp$ is given by
\[ 
\gp = \left\{A \in \gso_{r,r+k} \ \left| \  A = \begin{pmatrix} 0 & B \\ B^T & 0 \end{pmatrix},\ B \in M_{r,r+k}(\bbr) \right. \right\},
\]
where $M_{r,r+k}(\bbr)$ denotes the vector space of $r \times (r+k)$-matrices with real coefficients. We define a linear subspace $\gm$ of $\gp$ by
\[ 
\gm = \left\{ \left. A = \begin{pmatrix} 0 & B \\ B^T & 0 \end{pmatrix} \in \gp \ \right| \ B = \begin{pmatrix} C & 0  \end{pmatrix},\ C \in M_{r,r+k-1}(\bbr)  \right\}.
\]
A straightforward calculation shows that $[[\gm,\gm],\gm] \subset \gm$, that is, $\gm$ is a Lie triple system in $\gp$. We have $\dim(\gm) = r(r+k-1)$ and hence $\dim(\gp) - \dim(\gm) = r(r+k) - r(r+k-1) = r$. Thus the connected complete totally geodesic submanifold $\Sigma$ of $M$ corresponding to the Lie triple system ${\mathfrak m}$, which is isometric to $SO^o_{r,r+k-1}/SO_rSO_{r+k-1}$, satisfies $\codim(\Sigma) = r = \rk(M)$. From Theorem \ref{bound} it follows that the index of $SO^o_{r,r+k}/SO_rSO_{r+k}$ is equal to $r$.

Altogether we have now proved the ``if''-part of Theorem \ref{main}:

\begin{prop}\label{examples}
Let $M$ be one of the following Riemannian symmetric spaces of noncompact type:
\begin{itemize}[leftmargin=.3in]
\item[\rm{(i)}] $SL_{r+1}({\mathbb R})/SO_{r+1}$, $r \geq 1$;
\item[\rm{(ii)}] $SO^o_{r,r+k}/SO_rSO_{r+k}$, $r \geq 1$, $k \geq 0$, $(r,k) \notin \{(1,0),(2,0)\}$.
\end{itemize}
Then $\rk(M) = r = i(M)$.
\end{prop}

\section {The classification} \label{proof}

The following result was proved in \cite{BO} and will be used later.

\begin {teo}\label{flat}
Let $M$ be an irreducible Riemannian symmetric space, $\Sigma$ a connected totally geodesic submanifold of $M$ and $p \in \Sigma$.  Then there exists a maximal abelian subspace $\ga$ of $\gp$ such that $\ga$ is transversal to $T_p\Sigma$, that is, $\ga \cap T_p\Sigma = \{0\}$.
\end {teo}

Let $M = G/K$ be an irreducible Riemannian symmetric space of noncompact type and assume that $i(M) = r = \rk(M)$. Then  there exists a connected complete totally geodesic submanifold $\Sigma$ of $M$ with $p \in \Sigma$ such that $\codim(\Sigma) = r$. According to Theorem \ref{flat} there exists a maximal abelian subspace $\ga$ of $\gp$ such that $\ga$ is transversal to $T_p\Sigma$. Let $\Psi$ be the set of restricted roots with respect to $\ga$ and $\Lambda = \{\alpha_1,\ldots,\alpha_r\}$ be a set of simple roots for $\Psi$. The next result provides a necessary criterion for an irreducible Riemannian symmetric space $M$ with $\rk(M) \geq 2$ to satisfy the equality $\rk(M) = i(M)$. 

\begin{prop}\label{boundary_reduction} {\sc (Boundary Reduction)}
Let $M$ be an irreducible Riemannian symmetric space of noncompact type with $\rk(M) \geq 2$ and assume that the equality $\rk(M) = i(M)$ holds. Then every irreducible boundary component $B_\Phi$ of $M$ satisfies $\rk(B_\Phi) = i(B_\Phi)$.
\end{prop}

\begin{proof}
Let $\Sigma_\Phi$ be the connected complete totally geodesic submanifold of $F_\Phi$ corresponding to the Lie triple system $T_p\Sigma \cap T_pF_\Phi$. Since $T_pM = T_p\Sigma \oplus \ga$ (direct sum) and $\ga \subset T_pF_\Phi$, we have $T_pF_\Phi = T_p\Sigma_\Phi \oplus \ga$ (direct sum). Thus the codimension of $\Sigma_\Phi$ in $F_\Phi$ is equal to $\dim{\ga} = r = \rk(M)= \rk(F_\Phi)$. 

The orthogonal projection $(T_p\Sigma_\Phi)_{T_pB_\Phi}$ of the Lie triple system $T_p\Sigma_\Phi$ onto $T_pB_\Phi$ is a Lie triple system. Let $\Sigma'_\Phi$ be the connected complete totally geodesic submanifold of $B_\Phi$ corresponding to the Lie triple system $(T_p\Sigma_\Phi)_{T_pB_\Phi}$. Since $T_pF_\Phi = T_p\Sigma_\Phi \oplus \ga = T_pB_\Phi \oplus \ga_\Phi$ (direct sum) and $\ga = \ga^\Phi \oplus \ga_\Phi$, we have $T_pB_\Phi = T_p\Sigma'_\Phi \oplus \ga^\Phi$ (direct sum), which implies that the codimension of $\Sigma'_\Phi$ in $B_\Phi$ is equal to $\dim(\ga^\Phi) = \dim(\ga) - \dim(\ga_\Phi) = r - (r-|\Phi|) = |\Phi| = \rk(B_\Phi)$. This implies $i(B_\Phi) \leq \rk(B_\Phi)$. However, since $B_\Phi$ is irreducible, we also have $\rk(B_\Phi) \leq i(B_\Phi)$ by Theorem \ref{bound}. Altogether this implies $\rk(B_\Phi) = i(B_\Phi)$.
\end{proof}

We recall the following result from \cite{BO}:

\begin{teo}\label{classification} {\sc (Symmetric spaces with index $\leq 3$)}
 Let $M$ be an irreducible Riemannian symmetric space of noncompact type.
\begin{itemize}[leftmargin=.3in]
\item[\rm(1)] $i(M) = 1$ if and only if $M$ is isometric to 
\begin{itemize}
\item[\rm(i)] the real hyperbolic space ${\mathbb R}H^{k+1} = SO^o_{1,1+k}/SO_{1+k}$, $k \geq 1$.
\end{itemize}
\item[\rm(2)] $i(M) = 2$ if and only if $M$ is isometric to one of the following spaces:
\begin{itemize}
\item[\rm(i)] the complex hyperbolic space ${\mathbb C}H^{k+1} = SU_{1,1+k}/S(U_1U_{1+k})$, $k \geq 1$;
\item[\rm(ii)] the Grassmannian $SO^o_{2,2+k}/SO_2SO_{2+k}$, $k \geq 1$;
\item[\rm(iii)] the symmetric space $SL_3({\mathbb R})/SO_3$.
\end{itemize}
\item[\rm(3)] $i(M) = 3$ if and only if $M$ is isometric to one of the following spaces:
\begin{itemize}
\item[\rm(i)] the Grassmannian $SO^o_{3,3+k}/SO_3SO_{3+k}$, $k \geq 1$;
\item[\rm(ii)] the symmetric space $G^2_2/SO_4$; 
\item[\rm(iii)] the symmetric space $SL_3({\mathbb C})/SU_3$;
\item[\rm(iv)] the symmetric space $SL_4({\mathbb R})/SO_4$.
\end{itemize}
\end{itemize}
\end{teo}

The Riemannian symmetric spaces of noncompact type with $\rk(M) = 1 = i(M)$ are precisely the real hyperbolic spaces $SO^o_{1,1+k}/SO_{1+k}$, $k \geq 1$. The irreducible Riemannian symmetric spaces of noncompact type with $\rk(M) \geq 2$ whose rank one boundary components are all real hyperbolic spaces are precisely those for which the restricted root system is reduced, that is, is not of type ($BC_r$). From Proposition \ref{boundary_reduction} we therefore obtain:

\begin{cor} {\sc (Rank One Boundary Reduction)} \label{rank_one_reduction}
Let $M$ be an irreducible Riemannian symmetric space of noncompact type with $\rk(M) \geq 2$ and assume that $\rk(M) = i(M)$. Then the restricted root system of $M$ is not of type ($BC_r$).
\end{cor}

According to Theorem \ref{classification}, the Riemannian symmetric spaces of noncompact type with $\rk(M) = 2 = i(M)$ are precisely $SO^o_{2,2+k}/SO_2SO_{2+k}$, $k \geq 1$, and $SL_3({\mathbb R})/SO_3$. The corresponding Dynkin diagrams with multiplicities are
\[
\xy
\POS (0,0) *\cir<2pt>{} ="a",
(10,0) *\cir<2pt>{}="b",
(0,-3) *{1},
(10,-3) *{k},
\ar @2{->} "a";"b"
\endxy
\ \ \textrm{and}\ \
\xy
\POS (0,0) *\cir<2pt>{} ="a",
(10,0) *\cir<2pt>{}="b",
(0,-3) *{1},
(10,-3) *{1},
\ar @{-} "a";"b"
\endxy\ \ .
\]
We can easily extract from Table \ref{dynkin} the Dynkin diagrams of rank $\geq 3$ with multiplicities for which every connected subdiagram of rank $2$ is one of the above: 
\[
\xy
\POS (0,0) *\cir<2pt>{} ="a",
(7,0) *\cir<2pt>{}="b",
(14,0) *\cir<2pt>{}="c",
(21,0) *\cir<2pt>{}="d",
(0,-3) *{1},
(7,-3) *{1},
(14,-3) *{1},
(21,-3) *{1},
\ar @{-} "a";"b",
\ar @{.} "b";"c",
\ar @{-} "c";"d"
\endxy
\ \ ,\ \ 
\xy
\POS (0,0) *\cir<2pt>{} ="a",
(7,0) *\cir<2pt>{}="b",
(14,0) *\cir<2pt>{}="c",
(21,0) *\cir<2pt>{}="d",
(28,0) *\cir<2pt>{}="e",
(0,-3) *{1},
(7,-3) *{1},
(14,-3) *{1},
(21,-3) *{1},
(28,-3) *{k},
\ar @{-} "a";"b",
\ar @{.} "b";"c",
\ar @{-} "c";"d",
\ar @2{->} "d";"e"
\endxy
\ \ ,\ \ 
\xy
\POS (0,0) *\cir<2pt>{} ="a",
(7,0) *\cir<2pt>{}="b",
(14,0) *\cir<2pt>{}="c",
(21,0) *\cir<2pt>{}="d",
(28,2) *\cir<2pt>{}="e",
(28,-2) *\cir<2pt>{}="f",
(0,-3) *{1},
(7,-3) *{1},
(14,-3) *{1},
(21,-3) *{1},
(31,3) *{1},
(31,-3) *{1},
\ar @{-} "a";"b",
\ar @{.} "b";"c",
\ar @{-} "c";"d",
\ar @{-} "d";"e",
\ar @{-} "d";"f"
\endxy\ \ ,
\]
\[
\xy
\POS (0,0) *\cir<2pt>{} ="a",
(14,5) *\cir<2pt>{} = "b",
(7,0) *\cir<2pt>{}="c",
(14,0) *\cir<2pt>{}="d",
(21,0) *\cir<2pt>{}="e",
(28,0) *\cir<2pt>{}="f",
(0,-3) *{1},
(12,5) *{1},
(7,-3) *{1},
(14,-3) *{1},
(21,-3) *{1},
(28,-3) *{1},
\ar @{-} "a";"c",
\ar @{-} "c";"d",
\ar @{-} "b";"d",
\ar @{-} "d";"e",
\ar @{-} "e";"f",
\endxy
\ \ ,\ \ 
\xy
\POS (0,0) *\cir<2pt>{} ="a",
(14,5) *\cir<2pt>{} = "b",
(7,0) *\cir<2pt>{}="c",
(14,0) *\cir<2pt>{}="d",
(21,0) *\cir<2pt>{}="e",
(28,0) *\cir<2pt>{}="f",
(35,0) *\cir<2pt>{}="g",
(0,-3) *{1},
(12,5) *{1},
(7,-3) *{1},
(14,-3) *{1},
(21,-3) *{1},
(28,-3) *{1},
(35,-3) *{1},
\ar @{-} "a";"c",
\ar @{-} "c";"d",
\ar @{-} "b";"d",
\ar @{-} "d";"e",
\ar @{-} "e";"f",
\ar @{-} "f";"g",
\endxy\ \ ,
\]
\[
\xy
\POS (0,0) *\cir<2pt>{} ="a",
(14,5) *\cir<2pt>{} = "b",
(7,0) *\cir<2pt>{}="c",
(14,0) *\cir<2pt>{}="d",
(21,0) *\cir<2pt>{}="e",
(28,0) *\cir<2pt>{}="f",
(35,0) *\cir<2pt>{}="g",
(42,0) *\cir<2pt>{}="h",
(0,-3) *{1},
(12,5) *{1},
(7,-3) *{1},
(14,-3) *{1},
(21,-3) *{1},
(28,-3) *{1},
(35,-3) *{1},
(42,-3) *{1},
\ar @{-} "a";"c",
\ar @{-} "c";"d",
\ar @{-} "b";"d",
\ar @{-} "d";"e",
\ar @{-} "e";"f",
\ar @{-} "f";"g",
\ar @{-} "g";"h"
\endxy
\ \ ,\ \ 
\xy
\POS (0,0) *\cir<2pt>{} ="a",
(7,0) *\cir<2pt>{}="b",
(14,0) *\cir<2pt>{}="c",
(21,0) *\cir<2pt>{}="d",
(0,-3) *{1},
(7,-3) *{1},
(14,-3) *{1},
(21,-3) *{1},
\ar @{-} "a";"b",
\ar @2{->} "b";"c",
\ar @{-} "c";"d"
\endxy
\ \ .
\]

\medskip
From Proposition \ref{boundary_reduction} we thus obtain:

\begin{cor} {\sc (Rank Two Boundary Reduction)} \label{rank_two_reduction}
Let $M$ be an irreducible Riemannian symmetric space of noncompact type with $\rk(M) \geq 3$ and assume that $\rk(M) = i(M)$. Then $M$ must be among the following spaces:
\begin{itemize}[leftmargin=.3in]
\item[\rm(1)] $SL_{r+1}({\bbr})/SO_{r+1}$, $r \geq 3$;
\item[\rm(2)] $SO^o_{r,r+k}/SO_{r}SO_{r+k}$, $r \geq 3$, $k \geq 0$;
\item[\rm(3)] $E_6^6/Sp_4$;
\item[\rm(4)] $E_7^7/SU_8$;
\item[\rm(5)] $E_8^8/SO_{16}$;
\item[\rm(6)] $F_4^4/Sp_3Sp_1$.
\end{itemize}
\end{cor}

We know from Proposition \ref{examples} that the symmetric spaces in (1) and (2) satisfy the equality $\rk(M) = i(M)$. 
In order to prove Theorem \ref{main} it remains to show that the four exceptional spaces in Corollary \ref{rank_two_reduction} do not satisfy the equality $\rk(M) = i(M)$. For $M = F_4^4/Sp_3Sp_1$ we can apply rank three boundary reduction:

\begin{cor} \label{f4sp3sp1}
The symmetric space $M =  F_4^4/Sp_3Sp_1$ does not satisfy the equality $\rk(M) = i(M)$.
\end{cor}

\begin{proof}
The Dynkin diagram with multiplicities for $F_4^4/Sp_3Sp_1$ is
\[
\xy
\POS (0,0) *\cir<2pt>{} ="a",
(7,0) *\cir<2pt>{}="b",
(14,0) *\cir<2pt>{}="c",
(21,0) *\cir<2pt>{}="d",
(0,-3) *{1},
(7,-3) *{1},
(14,-3) *{1},
(21,-3) *{1},
\ar @{-} "a";"b",
\ar @2{->} "b";"c",
\ar @{-} "c";"d"
\endxy\ .
\]
We see from Theorem \ref{classification} that the boundary component $B_\Phi = Sp_3({\bbr})/U_3$ corresponding to the rank three subdiagram 
\[
\xy
\POS 
(7,0) *\cir<2pt>{}="b",
(14,0) *\cir<2pt>{}="c",
(21,0) *\cir<2pt>{}="d",
(7,-3) *{1},
(14,-3) *{1},
(21,-3) *{1},
\ar @2{->} "b";"c",
\ar @{-} "c";"d"
\endxy
\]
does not satisfy the equality $\rk(B_\Phi) = i(B_\Phi)$. The statement thus follows from Proposition \ref{boundary_reduction}.
\end{proof}

The situation for the exceptional symmetric space $E_6^6/Sp_4$ is quite interesting as the following result shows.

\begin{prop} \label{e6sp4}
Every irreducible boundary component $B_\Phi$ of $M = E_6^6/Sp_4$ satisfies $\rk(B_\Phi) = i(B_\Phi)$. However, $M$ does not satify the equality $\rk(M) = i(M)$.
\end{prop}

\begin{proof}
We list the different types of irreducible boundary components of $M$ by cardinality of $\Phi$.
\begin{itemize}[leftmargin=.3in]
\item[\rm(1)] $|\Phi| = 1$: $B_\Phi = SL_2(\bbr)/SO_2$;
\item[\rm(2)] $|\Phi| = 2$: $B_\Phi = SL_3(\bbr)/SO_3$;
\item[\rm(3)] $|\Phi| = 3$: $B_\Phi = SL_4(\bbr)/SO_4$;
\item[\rm(4)] $|\Phi| = 4$: $B_\Phi = SL_5(\bbr)/SO_5$ and $B_\Phi = SO^o_{4,4}/SO_4SO_4$;
\item[\rm(5)] $|\Phi| = 5$: $B_\Phi = SL_6(\bbr)/SO_6$ and $B_\Phi = SO^o_{5,5}/SO_5SO_5$.
\end{itemize}
As we have shown in Proposition \ref{examples}, each of these boundary components satisfies $\rk(B_\Phi) = i(B_\Phi)$.

We have $n = \dim(M) = 42$ and $r = \rk(M) = 6$. Assume that there exists a maximal totally geodesic submanifold $\Sigma$ of $M$ with $d = \codim(\Sigma) = 6$. We first assume that $\Sigma$ is semisimple. Then the inequality in Corollary \ref{ref4} is satisfied and thus $\Sigma$ is a reflective submanifold of $M$. As usual, we write $\Sigma = G'/K'$, where $G'$ is the connected closed subgroup of $E_6^6$ with Lie algebra $\gg' = [T_p\Sigma , T_p\Sigma] \oplus T_p\Sigma$ and $K' = G'_p$. Note that $K'$ is connected since $\Sigma$ is simply connected. Let $s \in I(M)$ be the geodesic reflection of $M$ in $\Sigma$ and define $\tau : E_6^6 \to E_6^6,\ g \mapsto sgs^{-1}$. It is clear that $G'$, and hence also $K'$, is contained in the fixed point set of $\tau$. Since $s$ commutes with the geodesic symmetry of $M$ at $p$, we have $\tau(Sp_4) = Sp_4$. Let $H$ be the connected component of the fixed point set of $\tau_{|Sp_4}$ containing the identity transformation of $Sp_4$. Note that $K' \subset H$. Then $Sp_4/H$ is a (simply connected) Riemannian symmetric space of compact type. However, as we observed in the proof of Corollary \ref {ref4}, we have $\dim(K ') \geq \dim(\Sigma) - \rk(M) = 30$ and therefore $\dim(Sp_4/H) \leq \dim(Sp_4/K') \leq 6$. Since there is no Riemannian symmetric space of $Sp_4$ of dimension $\leq 6$ we conclude that there is no reflective submanifold $\Sigma$ of $M$ with $\codim(\Sigma) = 6$. [Note: This fact can also be seen directly from Leung's classification of reflective submanifolds. However, we prefer to give a conceptual proof here.] Therefore $\Sigma$ cannot be semisimple. If $\Sigma$ is non-semisimple, then $\Sigma = \bbr \times SO^o_{5,5}/SO_5SO_5$ by Table \ref{totgeodnsstable} and hence $\codim(\Sigma) = 16$, which is a contradiction. Altogether we can now conclude that there is no totally geodesic submanifold in $M$ with $\codim(M) = 6$. This implies $rk(M) < i(M)$.
\end{proof}

As a consequence of Proposition \ref{e6sp4} we can now settle the two remaining cases.

\begin{cor} \label{e7e8}
The symmetric spaces $M =  E_7^7/SU_8$ and $M = E^8_8/SO_{16}$ do not satisfy the equality $\rk(M) = i(M)$.
\end{cor}

\begin{proof}
We see from Table \ref{dynkin} that the Dynkin diagram with multiplicities for $E_6^6/Sp_4$ can be embedded into the Dynkin diagrams with multiplicities for  $E_7^7/SU_8$ and $E^8_8/SO_{16}$. This means that $E_6^6/Sp_4$ is an irreducible boundary component of both $E_7^7/SU_8$ and $E^8_8/SO_{16}$. From Proposition \ref{boundary_reduction} and Proposition \ref{e6sp4} we can conclude that both  $M = E_7^7/SU_8$ and $M = E^8_8/SO_{16}$ do not satisfy the equality $\rk(M) = i(M)$.
\end{proof}

Theorem \ref{main} now follows from Proposition \ref{examples}, Corollary \ref{rank_two_reduction}, Corollary \ref{f4sp3sp1}, Proposition \ref{e6sp4} and Corollary \ref{e7e8}. We also obtain the following interesting characterization of the exceptional symmetric space $E_6^6/Sp_4$:

\begin{prop}
The exceptional symmetric space $E^6_6/Sp_4$ is the only irreducible Riemannian symmetric space $M$ of noncompact type with $rk(M) \geq 3$ for which every irreducible boundary component $B$ satisfies $\rk(B) = i(B)$ but the manifold itself does not satisfy $\rk(M) = i(M)$.
\end{prop}

\section{Further applications} \label{applications}

In this section we will calculate $i(M)$ for a few irreducible Riemannian symmetric spaces $M$ of noncompact type using the methods we developed in this paper and Leung's classification of reflective submanifolds. We first recall some known results to put our results into context.

The totally geodesic submanifolds of Riemannian symmetric spaces $M$ of noncompact type with $\rk(M) = 1$ were classified by Wolf in \cite{Wo}. We use the following notations: $\bbr H^{k+1} = SO^o_{1,1+k}/SO_{1+k}$ denotes the $(k+1)$-dimensional real hyperbolic space, $\bbc H^{k+1} = SU_{1,1+k}/S(U_1U_{1+k})$ denotes the $(k+1)$-dimensional complex hyperbolic space, $\bbh H^{k+1} = Sp_{1,1+k}/Sp_1Sp_{1+k}$ denotes the $(k+1)$-dimensional quaternionic hyperbolic space, and $\bbo H^2 = F_4^{-20}/Spin_9$ denotes the Cayley hyperbolic plane. Here, $k \geq 1$. The totally geodesic submanifolds of irreducible Riemannian symmetric spaces $M$ of noncompact type with $\rk(M) = 2$ were classified by Klein in \cite{K1}, \cite{K2}, \cite{K3} and \cite{K4}. From Wolf's and Klein's classifications we obtain $i(M)$ for all irreducible Riemannian symmetric spaces $M$ of noncompact type with $\rk(M) \leq 2$. Some of the indices for $\rk(M) = 2$ were calculated by Onishchik in \cite{On}. We summarize all this in Table \ref{iMforrkMleq2}. 

\begin{table}[h]
\caption{The index $i(M)$ for irreducible Riemannian symmetric spaces $M$ of noncompact type with $\rk(M) \leq 2$ and totally geodesic submanifolds $\Sigma$ of $M$ with $\codim(\Sigma) = i(M)$} 
\label{iMforrkMleq2} 
{\footnotesize\begin{tabular}{ | p{2.7cm}  p{4cm}  p{1.2cm}  p{0.8cm}  p{1.5cm} | }
\hline \rule{0pt}{4mm}
\hspace{-1mm}$M$ & $\Sigma$ & $\dim(M)$ & $i(M)$ & Comments\\[1mm]
\hline \rule{0pt}{4mm}
\hspace{-2mm} $\bbr H^{k+1}$ & $\bbr H^k$ & $k+1$ & $1$ & $k \geq 1$\\
$\bbc H^{k+1}$ & $\bbc H^k$ (and $\bbr H^2$ for $k=1$) & $2(k+1)$ & $2$ & $k \geq 1$\\
$\bbh H^{k+1}$ & $\bbh H^k$ (and $\bbc H^2$ for $k=1$)& $4(k+1)$ & $4$ & $k \geq 1$\\
$\bbo H^2$ & $\bbo H^1$, $\bbh H^2$ & $16$ & $8$ &  \\[1mm]
\hline \rule{0pt}{4mm}
\hspace{-2mm} 
$SL_3(\bbr)/SO_3$ & $\bbr \times \bbr H^2$ & $5$ & $2$ & \\
$SO^o_{2,2+k}/SO_2SO_{2+k}$ & $SO^o_{2,1+k}/SO_2SO_{1+k}$ & $2(k+2)$ &  $2$& $k \geq 1$\\
$SL_3(\bbc)/SU_3$ & $SL_3(\bbr)/SO_3$ & $8$ & $3$ & \\
$G_2^2/SO_4$ & $SL_3(\bbr)/SO_3$ & $8$ & $3$ &  \\
$SO_5(\bbc)/SO_5$ & $SO_4(\bbc)/SO_4$, $SO^o_{2,3}/SO_2SO_3$ & $10$ & $4$ &\\
$SU_{2,2+k}/S(U_2U_{2+k})$ & $SU_{2,1+k}/S(U_2U_{1+k})$ & $4(k+2)$ & $4$ &  $k \geq 1$\\
$SU^*_6/Sp_3$ & $SL_3(\bbc)/SU_3$, $\bbh H^2$ & $14$ & $6$ & \\
$G_2(\bbc)/G_2$ & $G_2^2/SO_4$, $SL_3(\bbc)/SU_3$ & $14$ & $6$ &\\
$Sp_{2,2}/Sp_2Sp_2$ & $Sp_2(\bbc)/Sp_2$ &  $16$ & $6$ & \\
$SO^*_{10}/U_5$ &$SO^*_8/U_4$, $SU_{2,3}/S(U_2U_3)$ & $20$ & $8$ & \\
$Sp_{2,2+k}/Sp_2Sp_{2+k}$ & $Sp_{2,1+k}/Sp_2Sp_{1+k}$ & $8(k+2)$ & $8$ &  $k \geq 1$\\
$E_6^{-26}/F_4$ & $\bbo H^2$ & $26$ & $10$ & \\
$E_6^{-14}/Spin_{10}U_1$ & $SO^*_{10}/U_5$ & $32$ & $12$ &  \\[1mm]
\hline
\end{tabular}}
\end{table}

Let $M$ be a connected Riemannian manifold and denote by ${\mathcal S}_r$ the set of all connected reflective submanifolds $\Sigma$ of $M$ with $\dim(\Sigma) < \dim(M)$. The reflective index $i_r(M)$ of $M$ is defined by
\[
i_r(M) = \min\{ \dim(M) - \dim(\Sigma) \mid \Sigma \in {\mathcal S}_r\} = \min\{ \codim(\Sigma) \mid \Sigma \in {\mathcal S}_r\}.
\]
It is clear that $i(M) \leq i_r(M)$ and thus $i_r(M)$ is an upper bound for $i(M)$. Leung classified in \cite{L1} and \cite{L2} the reflective submanifolds of irreducible simply connected Riemannian symmetric spaces of compact type. Using duality this allows us to calculate $i_r(M)$ explicitly for all irreducible Riemannian symmetric spaces $M$ of noncompact type. We list the reflective indices $i_r(M)$ for all $M$ with $\rk(M) \geq 3$ in Table \ref{summary}.

\begin{table}[h]
\caption{The reflective index $i_r(M)$ for irreducible Riemannian symmetric spaces $M$ of noncompact type with $\rk(M) \geq 3$ and  reflective submanifolds $\Sigma$ of $M$ with $\codim(\Sigma) = i_r(M)$} 
\label{summary} 
{\footnotesize\begin{tabular}{ | p{2.7cm}  p{3.7cm}  p{1.2cm}  p{0.8cm}  p{1.6cm}  p{2.2cm} |}
\hline \rule{0pt}{4mm}
\hspace{-1mm}$M$ & $\Sigma$ & $\dim(M)$ & $i_r(M)$ & Comments & $i(M) = i_r(M)$?\\[1mm]
\hline \rule{0pt}{4mm}
\hspace{-2mm} 
$SL_{r+1}(\bbr)/SO_{r+1}$ & $\bbr \times SL_r(\bbr)/SO_r$ & $\frac{1}{2}r(r+3)$ & $r$ & $r \geq 3$ & yes\\
$SL_4(\bbc)/SU_4$ & $Sp_2(\bbc)/Sp_2$ & $15$ & $5$ & & yes\\
$SL_{r+1}(\bbc)/SU_{r+1}$ & $\bbr \times SL_r(\bbc)/SU_r$ & $r(r+2)$ & $2r$ & $r \geq 4$ & ?\\
$SU^*_{2r+2}/Sp_{r+1}$ & $\bbr \times SU^*_{2r}/Sp_r$ & $r(2r+3)$ & $4r$ & $r \geq 3$ & ?\\[1mm]
\hline \rule{0pt}{4mm}
\hspace{-2mm} $SO^o_{r,r+k}/SO_{r}SO_{r+k}$ & $SO^o_{r,r+k-1}/SO_{r}SO_{r+k-1}$ & $r(r+k)$ &  $r$& $r \geq 3, k \geq 1$ & yes\\
$SO_{2r+1}(\bbc)/SO_{2r+1}$ & $SO_{2r}(\bbc)/SO_{2r}$ & $r(2r+1)$ & $2r$ & $r \geq 3$ & yes\\[1mm]
\hline \rule{0pt}{4mm}
\hspace{-2mm} 
$Sp_r(\bbr)/U_r$ & $\bbr H^2 \times Sp_{r-1}(\bbr)/U_{r-1}$ & $r(r+1)$ & $2r-2$ & $r \geq 3$ & yes for $r\leq 5$, otherwise ?\\
$SU_{r,r}/S(U_rU_r)$ & $SU_{r-1,r}/S(U_{r-1}U_r)$ & $2r^2$ & $2r$ & $r \geq 3$ & yes\\
$Sp_r(\bbc)/Sp_r$ & $\bbr H^3 \times Sp_{r-1}(\bbc)/Sp_{r-1}$ & $r(2r+1)$ & $4r-4$ & $r \geq 3$ & ?\\
$SO^*_{4r}/U_{2r}$ & $SO^*_{4r-2}/U_{2r-1}$ & $2r(2r-1)$ & $4r-2$ &  $r \geq 3$ & ?\\
$Sp_{r,r}/Sp_rSp_r$ & $Sp_{r-1,r}/Sp_{r-1}Sp_r$ &  $4r^2$ & $4r$ & $r \geq 3$ & ?\\
$E_7^{-25}/E_6U_1$ & $E_6^{-14}/Spin_{10}U_1$ & $54$ & $22$ & & ? \\[1mm]
\hline \rule{0pt}{4mm}
\hspace{-2mm} $SO^o_{r,r}/SO_{r}SO_{r}$ & $SO^o_{r-1,r}/SO_{r-1}SO_{r}$ & $r^2$ & $r$ &  $r \geq 4$ & yes\\
$SO_{2r}(\bbc)/SO_{2r}$ & $SO_{2r-1}(\bbc)/SO_{2r-1}$ & $r(2r-1)$ & $2r-1$ &  $r \geq 4$ & yes\\[1mm]
\hline \rule{0pt}{4mm}
\hspace{-2mm} $SU_{r,r+k}/S(U_rU_{r+k})$ & $SU_{r,r+k-1}/S(U_rU_{r+k-1})$ & $2r(r+k)$ & $2r$ &  $r \geq 3, k \geq 1$ & yes\\
$Sp_{r,r+k}/Sp_rSp_{r+k}$ & $Sp_{r,r+k-1}/Sp_rSp_{r+k-1}$ & $4r(r+k)$ & $4r$ &  $r \geq 3, k \geq 1$ & yes for $r-1 \leq k$, otherwise ?\\
$SO^*_{4r+2}/U_{2r+1}$ &$SO^*_{4r}/U_{2r}$ & $2r(2r+1)$ & $4r$ & $r \geq 3$ & ? \\[1mm]
\hline \rule{0pt}{4mm}
\hspace{-2mm}  $E_6^6/Sp_4$ & $F_4^4/Sp_3Sp_1$ &  $42$ & $14$ & & ?\\
$E_6(\bbc)/E_6$ & $F_4(\bbc)/F_4$ & $78$ &  $26$ & & ?\\[1mm]
\hline \rule{0pt}{4mm}
\hspace{-2mm} $E_7^7/SU_8$ & $\bbr \times E^6_6/Sp_4$ & $70$ & $27$ & & ?\\
$E_7({\mathbb C})/E_7$ & $\bbr \times E_6(\bbc)/E_6$ & $133$ & $54$ & & ?\\[1mm]
\hline \rule{0pt}{4mm}
\hspace{-2mm} $E_8^8/SO_{16}$ & $\bbr H^2 \times E_7^7/SU_8$ & $128$ & $56$ & & ?\\
$E_8(\bbc)/E_8$ & $\bbr H^3 \times E_7(\bbc)/E_7$ & $248$ & $112$ & & ?\\[1mm]
\hline \rule{0pt}{4mm}
\hspace{-2mm} $F_4^4/Sp_3Sp_1$ & $SO^o_{4,5}/SO_4SO_5$ & $28$ & $8$ & & yes\\
$E_6^2/SU_6Sp_1$ & $F_4^4/Sp_3Sp_1$ & $40$ & $12$ & & ? \\
$E_7^{-5}/SO_{12}Sp_1$ & $E_6^2/SU_6Sp_1$ & $64$ & $24$ & & ?\\
$E_8^{-24}/E_7Sp_1$ & $E_7^{-5}/SO_{12}Sp_1$ & $112$ & $48$ & & ?\\
$F_4(\bbc)/F_4$ & $SO_9(\bbc)/SO_9$ & $52$ & $16$&  & ?\\[1mm]
\hline
\end{tabular}}
\end{table}

As an application of Corollaries \ref{ref4} and \ref{nssrefl} we will now calculate the index of a few symmetric spaces. Let $\Sigma$ be a maximal totally geodesic submanifold of an $n$-dimensional irreducible Riemannian symmetric space $M$ of noncompact type with $r = \rk(M) \geq 2$ such that $i(M) = \codim(\Sigma)$. If $\Sigma$ is non-semisimple, then $\Sigma$ is a reflective submanifold by Corollary \ref{nssrefl}. If $\Sigma$ is semisimple and $d = \codim(\Sigma)$ satisfies $d(d+1) < 2(n - r)$, then $\Sigma$ is a reflective submanifold of $M$ by Corollary \ref{ref4}. It follows that if $\codim(\Sigma) \leq i_r(M) - 1$ and $(i_r(M)-1)i_r(M) < 2(n-r)$, then $\Sigma$ is a reflective submanifold. 
Altogether this implies the following

\begin{prop} \label{irM=iM}
Let $M$ be an irreducible Riemannian symmetric space of noncompact type with $\rk(M) \geq 2$. If
\[
(i_r(M)-1)i_r(M) < 2(\dim(M) - \rk(M)),
\]
then $i(M) = i_r(M)$.
\end{prop}

The inequality in Proposition \ref{irM=iM} can be checked explicitly for each symmetric space $M$ in Table \ref{summary}: 

\begin{cor}  \label{moreexamples}
The following Riemannian symmetric spaces $M$ of noncompact type with $\rk(M) \geq 3$ satisfy the inequality in Proposition \ref{irM=iM} and therefore satisfy the equality $i(M) = i_r(M)$:
\begin{itemize}[leftmargin=.3in]
\item[\rm(i)] $SL_{r+1}(\bbr)/SO_{r+1}$, $r \geq 3$;
\item[\rm(ii)] $SL_4(\bbc)/SU_4$;
\item[\rm(iii)] $SO^o_{r,r+k}/SO_rSO_{r+k}$, $r \geq 3$, $k \geq 1$;
\item[\rm(iv)] $SO_{2r+1}(\bbc)/SO_{2r+1}$, $r \geq 3$;
\item[\rm(v)] $Sp_r(\bbr)/U_r$, $3 \leq r \leq 4$;
\item[\rm(vi)] $SO^o_{r,r}/SO_rSO_r$, $r \geq 4$;
\item[\rm(vii)] $SO_{2r}(\bbc)/SO_{2r}$, $r \geq 4$;
\item[\rm(viii)] $SU_{r,r+k}/S(U_rU_{r+k})$, $r \geq 3$, $k \geq 1$;
\item[\rm(ix)] $Sp_{r,r+k}/Sp_rSp_{r+k}$, $3 \leq r \leq k$.
\end{itemize}
\end{cor}

We inserted this result into the last column of Table \ref{summary}.

We can also use these methods to determine all irreducible Riemannian symmetric spaces $M$ of noncompact type with $i(M) = 4$.

\begin{teo} \label{iM=4} {\sc (Symmetric spaces with index four)}
Let $M$ be an irreducible Riemannian symmetric space of noncompact type. Then $i(M) = 4$  if and only if $M$ is isometric to one of the following symmetric spaces:
\begin{itemize}[leftmargin=.3in]
\item[\rm{(i)}] $\bbh H^{k+1} = Sp_{1,1+k}/Sp_1Sp_k$, $k \geq 1$;
\item[\rm{(ii)}] $SU_{2,2+k}/S(U_2U_{2+k})$, $k \geq 1$;
\item[\rm{(iii)}] $SO^o_{4,4+k}/SO_4SO_{4+k}$, $k \geq 0$;
\item[\rm{(iv)}] $SO_5(\bbc)/SO_5$;
\item[\rm{(v)}] $Sp_3(\bbr)/U_3$;
\item[\rm{(vi)}] $SL_5(\bbr)/SO_5$.
\end{itemize}
\end{teo}

\begin{proof}
From Tables \ref{iMforrkMleq2} and \ref{summary} and Corollary \ref{moreexamples} we see that every symmetric space listed in Theorem \ref{iM=4} satisfies $i(M) = 4$. 
Conversely, let $M$ be an irreducible Riemannian symmetric space of noncompact type with $i(M) = 4$ and let $\Sigma$ be a maximal totally geodesic submanifold of $M$ with $d = \codim(\Sigma) = 4$. If $\rk(M) \leq 2$ we obtain from Table \ref{iMforrkMleq2} that $M$ is one of the spaces in (i), (ii) and (iv). Assume that $\rk(M) \geq 3$.
If $\Sigma$ is non-semisimple, then $\Sigma$ is reflective by Corollary \ref{nssrefl}. If $\Sigma$ is semisimple and $\dim(M) - \rk(M) \geq 11$, then $\Sigma$ is reflective by Corollary \ref{ref4}. Thus we have $i_r(M) = i(M) = 4$ if $\dim(M) - \rk(M) \geq 11$ and we can use Table \ref{summary} to see that $M$ is isometric to a space in (iii). The symmetric spaces $M$ with $\rk(M) \geq 3$ and $\dim(M) - \rk(M) < 11$ are $SL_4(\bbr)/SO_4$ and $SO^o_{3,4}/SO_3SO_4$ (which both have index $3$ by Theorem \ref{classification}), $Sp_3(\bbr)/U_3$ and $SL_5(\bbr)/SO_5$ (which both have index $4$ by Corollary \ref{moreexamples}). This concludes the proof of Theorem \ref{iM=4}
\end{proof}

The analogous argument does not work for index five. For example, $M = SU_{3,3}/S(U_3U_3)$ has $i_r(M) = 6$, but for $d= 5$ the inequality $d(d+1) < 2(\dim(M) - \rk(M))$ is not satisfied, so we can only conclude $i(M) \in \{5,6\}$ with our results so far. However, using the classification in \cite{BT} of cohomogeneity one actions on irreducible Riemannian symmetric spaces of noncompact type, we can improve the inequality in Corollary \ref{ref4} when $\codim(\Sigma) \geq 5$:

\begin{prop} \label{ref4plus}
Let $M$ be an $n$-dimensional irreducible Riemannian symmetric space of noncompact type with $r = \rk(M) \geq 2$ and let $\Sigma$ be a semisimple connected complete totally geodesic submanifold of $M$ with $\codim(\Sigma) = d \geq 5$. If \[d(d-1) < 2(n - r - 1),\]then $\Sigma$ is a reflective submanifold of $M$.
\end{prop}

\begin{proof}
As usual, we write $\Sigma = G'/K'$ and identify $SO_d$ with $SO(\nu_p\Sigma)$. Since $d \geq 5$ and any connected subgroup of $SO_d$ is totally geodesic in $SO_d$, we see from Corollary \ref{moreexamples} that the minimal codimension of a connected subgroup of $SO_d$ is equal to $d-1$, which is exactly the codimension of $SO_{d-1}$. A principal $K'$-orbit on $\Sigma$ has dimension $n-d-\rk(\Sigma)$, which implies $\dim(K') \geq n - d - \rk(\Sigma) \geq n - d - r$. Consequently, if $\frac{1}{2}(d-1)(d-2) < n - d - r$, then $\dim(K') > \frac{1}{2}(d-1)(d-2) = \dim(SO_{d-1}) $. The inequality $\frac{1}{2}(d-1)(d-2) < n - d - r$ is equivalent to $d(d-1) < 2(n - r - 1)$. If the kernel $\ker(\rho)$ of the slice representation $\rho : K' \to SO (\nu_p \Sigma)$ has positive dimension, then $\Sigma$ is a reflective submanifold by Proposition \ref{ref3}. If $\dim(\ker(\rho)) = 0$, then we must have $\gk' = \so_d$ and the action of $K'$ on the unit sphere in $\nu_p\Sigma$ is transitive. This implies that $\Sigma$ is a totally geodesic singular orbit of a cohomogeneity one action on $M$. It was proved in \cite{BT} that with five exceptions any such orbit is reflective. Three of the five exceptions do not satisfy the assumption $d \geq 5$. The remaining two exceptions are $\Sigma = G_2^\bbc/G_2$ in $M = SO_7(\bbc)/SO_7$ and $\Sigma = SL_3(\bbc)/SU_3$ in $M = G_2^\bbc/G_2$, and both do not satisfy the inequality  $d(d-1) < 2(n - r - 1)$. It follows that $\Sigma$ is reflective.
\end{proof}

Note that the assumption $d \geq 5$ in Proposition \ref{ref4plus} is essential. For example, $\Sigma = G_2^2/SO_4$ is a semisimple totally geodesic submanifold of $M = SO^o_{3,4}/SO_3SO_4$ with $d = \codim(\Sigma) = 4$. The inequality in Proposition \ref{ref4plus} is satisfied, but $\Sigma$ is non-reflective. For $d=3$ the totally geodesic submanifold $\Sigma = SL_3(\bbr)/SO_3$ in $M = G_2^2/SO_4$ provides a counterexample. 

From Proposition \ref{ref4plus} we obtain:

\begin{cor} \label{rankreflective}
Let $M$ be an irreducible Riemannian symmetric space of noncompact type and let $\Sigma$ be a semisimple connected complete totally geodesic submanifold of $M$ with $\codim(\Sigma) \geq 5$. If  $\codim(\Sigma) = \rk(M)$, then $\Sigma$ is a reflective submanifold of $M$.
\end{cor}

\begin{proof}
For $d = \codim(\Sigma) = \rk(M) = r$ the inequality in Proposition \ref{ref4plus} becomes $r^2 + r < 2n-2$. It is clear that $n = \dim(M) \geq \frac{1}{2}\#(R) + r$, where $(R)$ denotes the restricted root system of $M$. For every root system occuring here we have $r^2 + r \leq \#(R)$, with equality if and only if $(R) = (A_r)$. Altogether this implies $r^2 + r \leq \#(R) \leq 2n - 2r < 2n-2$ and hence $\Sigma$ is reflective by Proposition \ref{ref4plus}. 
\end{proof}

From Proposition \ref{ref4plus} we also obtain:

\begin{cor} \label{betterestimate}
Let $M$ be an irreducible Riemannian symmetric space of noncompact type with $\rk(M) \geq 2$, $i(M) \geq 5$ and $i_r(M) \geq 6$. If
\[
(i_r(M)-2)(i_r(M)-1) < 2(\dim(M) - \rk(M) - 1),
\]
then $i(M) = i_r(M)$.
\end{cor}

\begin{proof}
Let $\Sigma$ be a maximal totally geodesic submanifold of $M$ such that $d = \codim(\Sigma) = i(M) \geq 5$. We put $n = \dim(M)$ and $r = \rk(M)$. If $\Sigma$ is non-semisimple, then $\Sigma$ is a reflective submanifold by Corollary \ref{nssrefl} and hence $d \geq i_r(M)$. If $\Sigma$ is semisimple and $d < i_r(M)$, then $d(d-1) < 2(n - r - 1)$ by assumption and thus $\Sigma$ is a reflective submanifold by Corollary \ref{ref4plus}, which is a contradiction to $d < i_r(M)$. It follows that $d \geq i_r(M)$ and therefore $i(M) = i_r(M)$.
\end{proof}

We can use Corollary \ref{betterestimate}  to calculate a few more indices for symmetric spaces which cannot be obtained via the inequality in Proposition \ref{irM=iM} and are therefore not listed in Corollary \ref{moreexamples}:

\begin{cor} \label{moreindices}
The following symmetric spaces satisfy $i(M) = i_r(M)$:
\begin{itemize}[leftmargin=.3in]
\item[\rm(i)] $Sp_5(\bbr)/U_5$;
\item[\rm(ii)] $SU_{r,r}/S(U_rU_r)$, $r \geq 3$;
\item[\rm(iii)] $Sp_{r,r+k}/Sp_rSp_{r+k}$, $k + 1 = r \geq 3$;
\item[\rm(iv)] $F_4^4/Sp_3Sp_1$.
\end{itemize}
\end{cor}

\begin{proof}
Let $M$ be one of the symmetric spaces in (i)-(iv). It is clear that $rk(M) \geq 2$. From Theorems \ref{classification} and \ref{iM=4} we see that $i(M) \geq 5$ and from Table \ref{summary} we see that $i_r(M) \geq 6$. It is a straightforward calculation to show that $M$ satisfies the inequality in Corollary \ref{betterestimate}, which then implies $i(M) = i_r(M)$.
\end{proof}

We inserted this result into the last column of Table \ref{summary}.

We can now also settle the classifications for $i(M) = 5$ and $i(M) = 6$.

\begin{teo} \label{iM=5} {\sc (Symmetric spaces with index five)}
Let $M$ be an irreducible Riemannian symmetric space of noncompact type. Then $i(M) = 5$ if and only if $M$ is isometric to one of the following symmetric spaces:
\begin{itemize}[leftmargin=.3in]
\item[\rm{(i)}] $SO^o_{5,5+k}/SO_5SO_{5+k}$, $k \geq 0$;
\item[\rm{(ii)}] $SL_4(\bbc)/SU_4$;
\item[\rm{(iii)}] $SL_6(\bbr)/SO_6$.
\end{itemize}
\end{teo}

\begin{proof}
From Corollary \ref{moreexamples} and Table \ref{summary} we see that every symmetric space listed in Theorem \ref{iM=5} satisfies $i(M) = 5$. Conversely, let $M$ be an irreducible Riemannian symmetric space of noncompact type with $i(M) = 5$ and let $\Sigma$ be a maximal totally geodesic submanifold of $M$ with $d = \codim(\Sigma) = 5$. From Table \ref{iMforrkMleq2} we obtain $\rk(M) \geq 3$. If $\Sigma$ is non-semisimple, then $\Sigma$ is reflective by Corollary \ref{nssrefl}. If $\Sigma$ is semisimple and $\dim(M) - \rk(M) > 11$, then $\Sigma$ is reflective by Proposition \ref{ref4plus}. Thus we have $i_r(M) = i(M) = 5$ if $\dim(M) - \rk(M) > 11$ and we can use Table \ref{summary} to see that $M$ is isometric to one of the spaces in (i)-(iii). If $\dim(M) - \rk(M) < 11$ we saw in the proof of Theorem \ref{iM=4} that $i(M) \in \{3,4\}$. There is no symmetric space $M$ with $\rk(M) \geq 3$ and $\dim(M) - \rk(M) = 11$. This concludes the proof of Theorem \ref{iM=5}
\end{proof}

\begin{teo} \label{iM=6} {\sc (Symmetric spaces with index six)}
Let $M$ be an irreducible Riemannian symmetric space of noncompact type. Then $i(M) = 6$ if and only if $M$ is isometric to one of the following symmetric spaces:
\begin{itemize}[leftmargin=.3in]
\item[\rm{(i)}] $SO^o_{6,6+k}/SO_6SO_{6+k}$, $k \geq 0$;
\item[\rm{(ii)}] $SU_{3,3+k}/S(U_3U_{3+k})$, $k \geq 0$;
\item[\rm{(iii)}] $SU^*_6/Sp_3$;
\item[\rm{(iv)}] $G_2^\bbc/G_2$;
\item[\rm{(v)}] $Sp_{2,2}/Sp_2Sp_2$; 
\item[\rm{(vi)}] $Sp_4(\bbr)/U_4$; 
\item[\rm{(vii)}] $SO_7(\bbc)/SO_7$; 
\item[\rm{(viii)}] $SL_7(\bbr)/SO_7$.
\end{itemize}
\end{teo}

\begin{proof}
From Tables \ref{iMforrkMleq2} and \ref{summary} we see that every symmetric space listed in Theorem \ref{iM=6} satisfies $i(M) = 6$. Conversely, let $M$ be an irreducible Riemannian symmetric space of noncompact type with $i(M) = 6$ and let $\Sigma$ be a maximal totally geodesic submanifold of $M$ with $d = \codim(\Sigma) = 6$. If $\rk(M) \in \{1,2\}$ we see from Table \ref{iMforrkMleq2} that $M$ is one of the spaces in (iii)-(v).
We assume that $\rk(M) \geq 3$. If $\Sigma$ is non-semisimple, then $\Sigma$ is reflective by Corollary \ref{nssrefl}. If $\Sigma$ is semisimple and $\dim(M) - \rk(M) > 16$, then $\Sigma$ is reflective by Proposition \ref{ref4plus}. Thus we have $i_r(M) = i(M) = 6$ if $\dim(M) - \rk(M) > 16$ and we can use Table \ref{summary} to see that $M$ is isometric to one of the spaces in (i), (ii), (vii) and (viii). If $\dim(M) - \rk(M) < 12$ we saw in the proof of Theorem \ref{iM=5} that $i(M) \in \{3,4\}$. The symmetric spaces $M$ with $\rk(M) \geq 3$ and $12 \leq \dim(M) - rk(M) \leq 16$ are $SO^o_{3,5}/SO_3SO_5$ and $SO^o_{3,6}/SO_3SO_6$ (which both have index $3$ by Theorem \ref{classification}), $SO^o_{4,4}/SO_4SO_4$ and $SO^o_{4,5}/SO_4SO_5$ (which both have index $4$ by Theorem \ref{iM=4}), $SL_6(\bbr)/SO_6$ and $SL_4(\bbc)/SU_4$ (which both have index $5$ by Theorem \ref{iM=5}),   $Sp_4(\bbr)/U_4$ (which has index $6$ by Corollary \ref{moreexamples} and Table \ref{summary}), $SU_{3,3}/S(U_3U_3)$ (which has index $6$ by Corollary \ref{moreindices} and Table \ref{summary}).
This concludes the proof of Theorem \ref{iM=5}
\end{proof}

We cannot continue beyond $i(M) = 6$ with our methods. For example, the symmetric space $M = Sp_3(\bbc)/Sp_3$ satisfies $\dim(M) = 21$ and $\rk(M) = 3$. Thus the inequality $d(d-1) < 2(\dim(M) - \rk(M) - 1) = 34$ in Proposition \ref{ref4plus} is satisfied if and only if $d \leq 6$. However, from Table \ref{summary} we know that $i_r(M) = 8$. Thus we must have $i(M) \in \{7,8\}$. We cannot exclude the possiblity $i(M) = 7$ here. 

It is worthwhile to point out that the only irreducible Riemannian symmetric space $M$ with $i(M) < i_r(M)$ known to us is $M = G_2^2/SO_4$. This leads us to the

\smallskip
{\sc Conjecture.} Let $M$ be an irreducible Riemannian symmetric space of noncompact type and $M \neq G_2^2/SO_4$. Then $i(M) = i_r(M)$.

\smallskip
We verified the conjecture in this paper for several symmetric spaces and for all symmetric spaces with $i(M) \leq 6$ or $\dim(M) \leq 20$. In the last column of Table \ref{summary} we summarize the current status of this conjecture.

\end {document}